\documentclass{article}
\usepackage[utf8]{inputenc}
\usepackage[a4paper,margin=2cm]{geometry}
\usepackage{graphicx,subcaption,amsmath,amssymb, comment,color,hyperref}
\graphicspath{{./Figures/}}
\usepackage{microtype}
\usepackage{enumitem}
\usepackage{amsthm}
\usepackage[final]{changes}

\newcommand{\gp}{Gaussian process}
\newcommand{\gps}{Gaussian processes}
\newcommand{\sa}{sensitivity analysis}
\newcommand{\si}{sensitivity index}
\newcommand{\sis}{sensitivity indices}

\title{Efficient estimation of divergence-based \sis{} with \gp{} surrogates}
\author{A.W. Eggels, D.T. Crommelin}

\begin{document}

\maketitle
\begin{abstract}	
	We consider the estimation of sensitivity indices based on divergence measures such as Hellinger distance. For sensitivity analysis of complex models, these divergence-based indices can be estimated by Monte-Carlo sampling (MCS) in combination with kernel density estimation (KDE). In a direct approach, the complex model must be evaluated at every input point generated by MCS, resulting in samples in the input-output space that can be used for density estimation. However, if the computational cost of the complex model strongly limits the number of model evaluations, this direct method gives large errors. We propose to use Gaussian process (GP) surrogates to increase the number of samples in the combined input-output space. By enlarging this sample set, the KDE becomes more accurate, leading to improved estimates. To compare the GP surrogates, we use a surrogate constructed by samples obtained with stochastic collocation, combined with Lagrange interpolation. Furthermore, we propose a new estimation method for these sensitivity indices based on minimum spanning trees. Finally, we also propose a new type of sensitivity indices based on divergence measures, namely direct sensitivity indices. These are useful when the input data is dependent.

\end{abstract}

%


\section{Introduction}
Sensitivity analysis is an essential part of uncertainty quantification and a very active research field \cite{Sal04,2Oak04,2Sal10}. Several types of sensitivity indices have been formulated, such as variance-based (including Sobol's indices \cite{Sob93}), density-based \cite{2Bor07}, derivative-based \cite{4Bli14} or divergence-based. Broadly speaking, divergence-based \sis{} quantify the difference between the joint probability distribution (or density) of model input and output on the one hand, and the product of their marginal distributions on the other hand.  A variety of divergence-based indices can be brought in a common framework built on the notion of $f$-divergence \cite{4Csi04}, as was shown by Da Veiga \cite{4Dav14}. The $f$-divergence is a generalization of several well-known divergences such as the Kullback-Leibler divergence \cite{4Erv14} and the Hellinger distance \cite{4Gib02}.

In most cases, these \sis{} cannot be computed analytically because the distribution of the model output given the input is not known exactly. As an alternative, one can resort to Monte Carlo (MC)  sampling combined with kernel density estimation: the input distribution is sampled using MC, the model is evaluated on all sampled input points, and from resulting input-output points the joint and marginal probability densities of input and output are estimated. However, when the number of available output points is low, for example because of high computational cost of the model, the estimated densities will generally be inaccurate, resulting in large errors in the estimated \sis{}.

In this study we first propose to increase the number of output samples by using a \gp{} (GP) surrogate. The GP is constructed on the input-output points that are obtained with the (expensive) model. The main idea is that the additional output samples improve the kernel density estimates even though they introduce a bias due to the difference between the true model and its GP approximation. Our approach is based on both the development of divergence-based indices and the use of \gps{} in \sa{}. Therefore, we briefly summarize some of the advancements in these areas. 
Auder \& Iooss \cite{2Aud08} presented two \sa{} methods based on Shannon and Kullback-Leiber entropy, respectively, building on work in \cite{2Krz01} and \cite{2Liu06}. Da Veiga \cite{4Dav14} introduced \sis{} based on the $f$-divergence. Recently, KDE also appears in estimators of mutual information measures in \cite{4Moo17}, where $f$-divergences are computed between the joint distribution of two random variables  and the product of their marginal distributions. In \cite{4Nos17}, $f$-divergence measures are computed by a $k$-nearest neighbor graph. 

The use of GPs is discussed in Marrel et al. \cite{4Mar09}, together with the analytical expressions for Sobol indices that arise from them. To compute the indices, two approaches are considered: one in which the predictor of the GP is used and one in which the full GP is used. The latter approach is found to be superior in convergence and robustness. Furthermore, the modeling error of the GP is integrated through confidence intervals; it is reported that the bias due to the use of the GP is negligible \cite{4Mar09}. In a related study, Svenson et al. \cite{4Sve14} estimate Sobol indices with GPs, using specific compactly supported kernel functions. Furthermore, combining GPs with derivative-based indices has been investigated by \cite{4Bli14} and \cite{4Loz16}. In \cite{4Paa17}, predictions from a GP are used to rank the input variables based on their predictive relevance. Two methods for this are presented in \cite{4Paa17}, one based on Kullback-Leibler divergence and one based on the variance of the posterior mean. 

Despite the developments sketched above, approaches that combine GP surrogate modeling and divergence-based \sa{} have not been explored much yet, although \cite{TBor12} already applied this approach. The methodology proposed in this paper combines these two elements. 

We note that for the approach proposed here it is not needed to assume that the inputs are mutually independent, nor does dependency of inputs make it more complicated. We present test cases with independent inputs as well as cases with dependent inputs. \added{For the former, we compare with results obtained with stochastic collocation \cite{2Mat03,Xiu05}. In this method, an appropriate set of points, called collocation points, is obtained. These are usually chosen as the zeros of the orthogonal polynomials with respect to the marginal input probability distributions. Then, Lagrange interpolation is used to approximate the output function. For dependent inputs, this method might not be ideal, as \cite{Nav15} already showed. }

\added{Second, we propose a new estimation method for the divergence-based \sis{} as introduced before. Because the KDE method depends on the choice of both kernel and kernel bandwidth, we propose to use an estimator without parameters which is numerically fast as well. This estimator is based on the approximation of one of the integrals appearing in the \si{} by computing a minimum spanning tree \cite{2Her02}. }

\added{As a third contribution, we propose a new set of \sis{} to complement the ones introduced before. This new set computes the direct \sis{}, which measure the sensitivity of the output with respect to one input variable only. This is beneficial for cases when the input variables are dependent, because these indices remove indirect effects caused by dependent input variables. To illustrate this, consider an example where $X=(X^1,X^2)$ follows a bivariate normal distribution with means $0$, variances $1$ and covariance $\rho>0$, while $u(x_1,x_2)=x_1$. Then the direct effect of $X^2$ on the output is zero, while the original \si{} would be positive due to the dependence between $X^1$ and $X^2$.}

Section \ref{sec:theory} describes the \sis{} central to this paper, their estimation method and the complications therein. It also contains our proposed method to enlarge the set of input and output data and the new estimation method. Section \ref{sec:results} applies these estimators to several test cases. Section \ref{sec:conclusion} concludes.

\section{Divergence-based \sis{} and their estimation}\label{sec:theory}
\added{We start by introducing the \sis{} derived from the $f$-divergences in Section \ref{ssec:sis}. Section \ref{ssec:problem} discusses the complications in estimating them. Gaussian processes and the two estimators are given in Section \ref{ssec:imp}. }

\subsection{Sensitivity indices from the \texorpdfstring{$f$}{f}-divergence}\label{ssec:sis}
We consider the situation where a model takes a vector of inputs $(X^1, ..., X^d)$ and returns a (scalar) output $Y$. The input vector $X$ is random, and as a result the output $Y$ is a random variable as well.
Da Veiga \cite{4Dav14} proposed to perform global sensitivity analysis with dependence measures, especially $f$-divergences (see also \cite{4Rah16}). In this way, the impact of the $k$th input variable $X^k$ on the output $Y$ is given by
\begin{equation}\label{eq:si}
	S_{X^k} = \mathbb{E}\left[d(Y,Y|X^k)\right],
\end{equation}
where $d(\cdot,\cdot)$ denotes a dissimilarity measure. The unnormalized first-order Sobol indices can also be written in this framework, namely with
\begin{equation*}
	d(Y,Y|X^k) = \left(\mathbb{E}(Y)-\mathbb{E}(Y|X^k)\right)^2.
\end{equation*}
We will use the Csisz\'ar $f$-divergence \cite{4Csi04}, which is given by
\begin{equation}\label{eq:f}
	d_f(Y,Y|X^k) = \int_\mathbb{R}f\left(\frac{p_Y(y)}{p_{Y|X^k}(y)}\right)p_{Y|X^k}(y) dy,
\end{equation}
with $f(\cdot)$ a convex function with $f(1)=0$, and $p_\cdot(\cdot)$ denotes a probability distribution function. Some well-known choices for $f$ are $f(t)=-\log(t)$ (Kullback-Leibler divergence) and $f(t)=(\sqrt{t}-1)^2$ (Hellinger distance). Combining (\ref{eq:si}) and (\ref{eq:f}) with basic probability theory gives us
\begin{equation}\label{eq:sif}
	S^f_{X^k} = \iint_{\mathbb{R}^2}f\left(\frac{p_Y(y)p_{X^k}(x)}{p_{X^k,Y}(x,y)}\right)p_{X^k,Y}(x,y) dydx.
\end{equation}
These sensitivity indices are equal to zero for $X^k$ and $Y$ independent and positive otherwise. Furthermore, they are invariant with respect to smooth and uniquely invertible transformation of $X^k$ and $Y$ \cite{2Kra04}, in contrast to Sobol indices which are only invariant with respect to linear transformations. Moreover, it is easy to generalize (\ref{eq:sif}) to multidimensional $X^{k,l}$. 

\subsection{Difficulties for estimation}\label{ssec:problem}
The main problem for computing $S^f_{X^k}$ is that the probability densities in (\ref{eq:sif}) are not known. In order to estimate $S^f_{X^k}$ it is necessary to estimate $p_Y(\cdot)$ and $p_{X^k,Y}(\cdot,\cdot)$, and, depending on the type of input, $p_{X^k}(\cdot)$ as well. In \cite{4Dav14} it is indicated that if samples $(X_L,Y_L)$ are available, only the ratio $r(x,y)=\frac{p_Y(y)p_{X^k}(x)}{p_{X^k,Y}(x,y)}$ needs to be estimated. 

The estimates of the densities can be obtained with kernel-density estimation (also in \cite{4Dav14,4Rah16}). To do so, one chooses a suitable kernel and a suitable value for the kernel bandwidth $h$. When the density of the input $X$ is known, this information can be used to determine $h$, otherwise, guidelines are available \cite{4Sco92}.

Clearly, the estimate of the density $p_Y$ will not be perfect, leading to an error in the estimation of $S^f_{X^k}$. This is strongly related to the number of samples $(X_L,Y_L)$ available for density estimation. If high computational cost of the model limits this number, the estimation of $S^f_{X^k}$ can be improved by using a surrogate of the model to generate more samples. \added{One possible way to do so is to use stochastic collocation (SC) \cite{2Mat03,Xiu05}. Herein, one chooses the samples $X_L$ as the collocation points, which are obtained as the zeros of the orthogonal polynomials with respect to the marginal input distributions. Then, at these collocation points, the corresponding output samples are obtained. Finally, an emulator is constructed by Lagrange interpolation on these samples.}

As an alternative, we propose to use Gaussian processes \cite{3Ras06} as a surrogate model to obtain the larger sample $(X_+,Y_+)=(X_L\cup X_{L^+},Y_L\cup Y_{L^+})$, in which $Y_{L^+}$ indicates the surrogate model output for the extra input samples $X_{L^+}$. For each data point in $X_{L^+}$, this $Y_{L^+}$ is a normal distribution in itself, and for each point in $X_L$ it is a degenerated normal distribution (i.e., it has zero variance). An additional advantage may be the availability of confidence intervals for $S^f_{X^k}$ at almost no extra computational cost. Unfortunately, these confidence intervals do not include the bias from approximating the output by a \gp{}. 

\subsection{Estimation using \gps{}}\label{ssec:imp}
We assume the input samples $X_L:=\{\mathbf{x}_l\}_{l=1}^L$ are already available, otherwise one can use Monte Carlo sampling (or Latin hypercube sampling in the case of independent uniform data) to select samples from the data $X$. Although it may be tempting to use other sample selection methods, it is not guaranteed that they represent the distribution just as naive sampling would. Then, the corresponding output $Y_L:=\{y_l\}_{l=1}^L$ can be obtained as $Y_L=G(X_L)$ with $G$ the process to generate output, which is either a function or a computational model. Then, one needs to fit a \gp{} $\widetilde{G}_{\{X_L,Y_L\}}(\mathbf{x})=N(\mu(\mathbf{x}),\Sigma(\mathbf{x}))$ to $(X_L,Y_L)$, thereby choosing an appropriate kernel. This \gp{} is now used to obtain output $Y_{L^+}=\widetilde{G}_{\{X_L,Y_L\}}(X_{L^+})$ for other input samples $X_{L^+}$. This leads to the augmented dataset $X_+=X_L\cup X_{L^+}$ of size $N=L+L_+$ with (partial) surrogate output $Y_+=Y_L\cup Y_{L^+}$. Note that $Y_{L^+}$ does not consist of single values, but rather of multivariate normal distributions. 

\subsubsection{Kernel density estimation}
We now explain how to compute the KDE on $(X_+,Y_+)$ and how it is used to approximate (\ref{eq:sif}). Because $S_{X^k}^f$ is computed per input variable $X^k$, it is here enough to consider one-dimensional kernel densities. 

For each input variable $X^k$ and output variable $Y$, the estimators for the kernel density are given by \cite{4Rah16}:
\begin{align*}
\widehat{f_{X^k}}(x) & = \frac{1}{Jh_{X^k}}\sum_{j=1}^J K_{X^k}\left(\frac{x-x_j}{h_{X^k}}\right), \\
\widehat{f_Y}(y) & = \frac{1}{Jh_Y}\sum_{j=1}^J K_Y\left(\frac{y-y_j}{h_Y}\right), \\
\widehat{f_{{X^k},Y}}(x,y) & = \frac{1}{Jh_{X^k}h_Y}\sum_{j=1}^J K_{X^k}\left(\frac{x-x_j}{h_{X^k}}\right)K_Y\left(\frac{y-y_j}{h_Y}\right),
\end{align*}
with $(x_j,y_j)$ the $j$th sample of the input data $(X^k,Y)$ and $J$ the size of the data. Note the input data $X=(X^1,\ldots,X^d)$ has to represent the distribution of $X$. An extension to a higher-dimensional $X^k$ is easy to obtain. For our purpose, we either have $J=L$ and $(X,Y)=(X_L,Y_L)$, or we have $J=N$ and $(X,Y)=(X_+,Y_+)$. We choose the Gaussian kernel and $h_{X^k}=h_Y=h$ according to Scott's rule \cite[p.\ 152]{TSco92}, which is optimized with respect to the normal distribution. Then, the estimator for $S^f_{X^k}$ as given by \cite{4Rah16} is obtained:
\begin{equation}\label{eq:H}
\overline{H}_{X^k,f}^{(J)} := \frac{1}{J}\sum_{j=1}^J f\left(\frac{\widehat{f_{X^k}}(x_j)\widehat{f_Y}(y_j)}{\widehat{f_{{X^k},Y}}(x_j,y_j)}\right).
\end{equation}
We note this choice of $h$ may not be optimal. We have adapted the bandwidth $h$ previously to the ranges of $X$ and $Y$, but the results of this are worse than with a single bandwidth. Also, kernel density estimation may not be the best choice when the domain of a variable $X^k$ or $Y$ is bounded and this variable has nonzero density at the boundaries.  

Until so far, we ignored the fact $Y_{L^+}$ is a multivariate normal random variable instead of a single value when $J=N$. Therefore, there are two options to obtain values for $Y_{L^+}$. The first option is to use the prediction mean $\mu(x)$ and get the resulting output samples 
\begin{equation}\label{eq:est1}
Y_{L^+}=\mu\left(X_{L^+}\right),
\end{equation}
to be used in (\ref{eq:H}). The other is to sample from this normal distribution $n_s$ times. In that case, one gets the $n_s$ output sets
\begin{equation}\label{eq:est2}
Y_{L^+}^{(s)}\sim N\left(\mu\left(X_{L^+}\right),\Sigma\left(X_{L^+}\right)\right),
\end{equation}
in which $\sim$ denotes ``sampled from the distribution'', and thereby $n_s$ estimates of $\overline{H}_{X^k,f}^{(N)}$. Note that this also implies the kernel density estimates have to be computed $n_s$ times. Because the computation of the kernel density estimate is expensive, we choose not to include this option. \added{We will indicate the estimator $\overline{H}_{X^k,f}^{(J)}$ by $\widehat{S^f_{X^k}}$ in the results, where the value of $J$ is clear from the context.}

\subsubsection{Minimum spanning trees}\label{ssec:cod}
\added{Before we can explain how to use the minimum spanning trees, we first need to introduce the concept of R\'enyi entropy. This is a generalization of the continuous Shannon entropy (see e.g. \cite{2Cov06}) and is given by}
\begin{equation*}
	H_\alpha(X) = \frac{1}{1-\alpha}\log\left(\int_\Omega (p(x))^\alpha \mathrm{d}x\right),
\end{equation*}
\added{for $\alpha\in(0,\infty)$. In the limit of $\alpha\rightarrow 1$, the R\'enyi entropy converges to the continuous Shannon entropy. Hero \& Michel \cite{2Her98,2Her02,2Her03} proposed a direct way to estimate the R\'enyi entropy for $\alpha\in(0,1)$ given a dataset $X_N$ consisting of $N$ samples of the probability distribution $X$ of dimension $d$. Their estimator is }
\begin{equation}\label{eq:Renyiest}
	\hat{H}_\alpha (X_N) = \frac{1}{1-\alpha} \left(\log \left(\frac{L_\gamma(X_N)}{N^\alpha}\right) - \log \beta_{L,\gamma}\right) = \frac{1}{1-\alpha} \log \left(\frac{L_\gamma(X_N)}{\beta_{L,\gamma} \, N^\alpha}\right),
\end{equation}
\added{in which $\gamma$ can be derived from the relation $\alpha=(d-\gamma)/d$. The functional $L_\gamma(X_N)$ is given by}
\begin{equation}\label{eq:lgamma}
L_\gamma(X_N) = \min_{T(X_N)} \sum_{e\in T(X_N)} |e|^\gamma,
\end{equation}
\added{in which $T(X_N)$ denotes the set of spanning trees on $X_N$ and $e$ denotes an edge therein. The parameter $\gamma$ can be computed from the desired value of $\alpha$ and will be within the interval $(0,d)$. The constant
$\beta_{L,\gamma}$ is defined by }
\begin{equation}\label{eq:beta}
	\beta_{L,\gamma} = \lim_{N\rightarrow \infty} L_\gamma(X_N)/N^\alpha,
\end{equation}
\added{for $X_N$ a sample of size $N$ of the uniform distribution in $d$ dimensions. However, we will estimate it for $N$ samples only by computing it for several repetitions of the sample $X_N$. }

\added{The estimator (\ref{eq:Renyiest}) is asymptotically unbiased and strongly consistent for $\alpha \in (0,1)$ \cite{2Her99b}. We focus on the case $\alpha=1/2$ and $d=2$ wherein $|e|$ denotes the Euclidean distance, hence $\gamma=1$. To see why we choose $\alpha=1/2$, we give the following derivation. First, we need to introduce the R\'enyi divergence by }
\begin{equation*}
	D_\alpha(f,g) = \frac{1}{\alpha-1}\log\left(\int_\Omega \left(\frac{f(x)}{g(x)}\right)^\alpha g(x)\mathrm{d}x\right),
\end{equation*}	
\added{for the probability distribution functions $f(\cdot)$ and $g(\cdot)$. We choose $f(\cdot)=p_{XY}(x,y)$ and $g(\cdot)=p_X(x)p_Y(y)$, where $p_{XY}$ is the joint probability distribution function of $X$ and $Y$ and $p_X(x)$ and $p_Y(y)$ denote the marginal probability distribution functions. Then, their R\'enyi divergence is given by}
\begin{align}
	D_{1/2}(p_{XY},p_Xp_Y) & = -2\log\left(\int_Y\int_X \left(\frac{p_{XY}(x,y)}{p_X(x)p_Y(y)}\right)^{1/2} p_X(x)p_Y(y)\mathrm{d}x\mathrm{d}y\right), \nonumber \\
	& = -2\log\left(\int_Y\int_X \sqrt{p_{XY}(x,y)p_X(x)p_Y(y)}\mathrm{d}x\mathrm{d}y\right) \label{eq:renyidiv}.
\end{align}
\added{We also have}
\begin{align}
	D_{1/2}(p_{XY},p_Xp_Y) & = -2\log\left(\int_Y\int_X \left(\frac{p_{XY}(x,y)}{p_X(x)p_Y(y)}\right)^{1/2} p_X(x)p_Y(y)\mathrm{d}x\mathrm{d}y\right), \nonumber \\
	& = -2\log\left(\int_Y\int_X \left(h_{xy}(x',y')\right)^{1/2} \mathrm{d}x'\mathrm{d}y'\right), \nonumber \\
	& = -H_{1/2}(h), \label{eq:dh}
\end{align}
\added{with $h(\cdot)$ the (well-defined) probability distribution function given by}
\begin{equation*}
	h(x,y) = \frac{p_{XY}(x,y)}{p_X(x)p_Y(y)}.
\end{equation*}

We also see that $S^H_{X^k}$, with $H$ denoting the \si{} derived from the Hellinger distance, is given through (\ref{eq:sif}) by
\begin{equation*}
S^H_{X^k} = \iint_{\mathbb{R}^2}\left(\sqrt{\frac{p_{X^k}(x)p_Y(y)}{p_{X^k,Y}(x,y)}}-1\right)^2p_{X^k,Y}(x,y) \mathrm{d}y\mathrm{d}x,
\end{equation*}
which can be simplified to
\begin{equation}\label{eq:sih}
S^H_{X^k} = 2 - 2 \iint_{\mathbb{R}^2}\sqrt{p_{X^k}(x)p_Y(y)p_{X^k,Y}(x,y)}\mathrm{d}y\mathrm{d}x.
\end{equation}
We now see the agreement between (\ref{eq:renyidiv}) and (\ref{eq:sih}). In case the domain of $X^k$ and $Y$ is extended to $\mathbb{R}$ by zero density outside of the domain, it is possible to write
\begin{equation*}
D_{1/2}(p_{X^kY},p_{X^k}p_Y) = -2\log(I), \quad S^H_{X^k} = 2-2I,
\end{equation*}
with 
\begin{equation*}
I = \iint_{\mathbb{R}^2}\sqrt{p_{X^k}(x)p_Y(y)p_{X^k,Y}(x,y)}\mathrm{d}y\mathrm{d}x,
\end{equation*}
hence 
\begin{equation*}
S^H_{X^k} = 2-2\exp\left(\frac{-D_{1/2}(p_{X^k,Y},p_{X^k}p_Y)}{2}\right). 
\end{equation*}
We can compute $S^H_{X^k}$ via $D_{1/2}(p_{X^k,Y},p_{X^k}p_Y)=-H_{1/2}(h)$ (Equation \ref{eq:dh}). Therefore, we need to estimate $L_\gamma(X^k,Y)$ (\ref{eq:lgamma}) and $\beta$ (\ref{eq:beta}). Because $I$ can be estimated as  
\begin{equation*}
\frac{L_\gamma(X^k,Y)}{\beta\sqrt{N}},
\end{equation*}
we estimate the \sis{} by
\begin{equation}\label{eq:SImst}
\widehat{S^H_{X^k}} = 2-2\frac{L_\gamma}{\beta\sqrt{N}}.
\end{equation}

\section{Results}\label{sec:results}
We test the estimators in several ways. The first test case is with regard to random input/output data and is described in Section \ref{ssec:random}. In this case, the estimates should be near zero. The second test case, in Section \ref{ssec:analytic}, is based on comparing analytic to numerical values of the \sis. In the third test case, the Ishigami function is used and tests are performed for both independent and dependent input data, of which the results can be found in Section \ref{ssec:Ishi}. The last test case is higher-dimensional and considers the Piston function (Section \ref{ssec:Piston}). In these tests, we only use the Hellinger distance. All experiments have been performed $n_r=10$ times. The error bars in the upcoming figures indicate the minimum and maximum value found. The results are summarized in Section \ref{ssec:rec}.

\subsection{Random data}\label{ssec:random}
\begin{figure}
	\centering
	\includegraphics[width=0.6\textwidth]{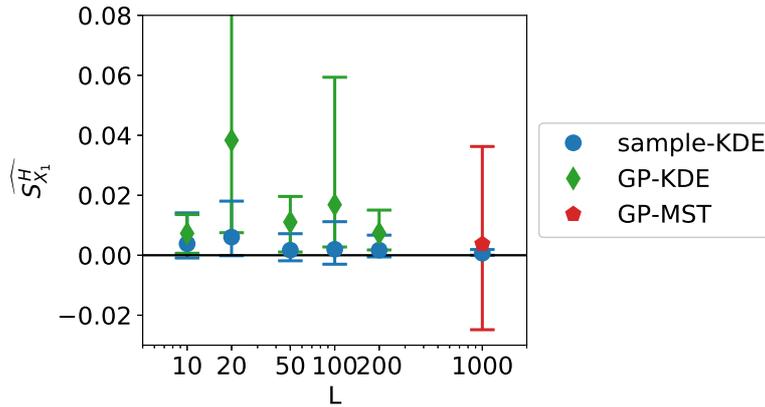}
	\caption{Sensitivity indices for random output.}\label{fig:zero}
\end{figure}

First, we check the behavior for random output, in which case the \sis{} should be zero. Both the input and output data are one-dimensional, uniformly distributed on $[0,1]$ and have size $N=10^3$, while $L$ is varied from $L=10$ to $L=200$ based on \cite{3Loe09}. The results are in Figure \ref{fig:zero}. 
On the right, we show the \si{} as computed on the complete, i.e., $L=N$, data by KDE and MST (blue circle and red pentagon). Herein, no \gp{} is used. As expected, their mean is around zero. The spread for the KDE method is smaller than for the MST method. The estimates for $\widehat{S^H_{X^k}}$ based on $L$ samples (blue circles) are also around zero, although their spread is larger than for $L=N$. Note that due to the numerical implementation, the \sis{} can become negative. 

We continue with the estimates based on \gps{}. Herein, the situation is a little different because the \gp{} fits a function through the data while there is no functional relation between input and output. Hence, the \sis{} will most likely not be equal to zero. When fitting the \gp{}, two cases appear, which have the same effect. The length scale and the process variance are either both small or both large. As a result, the predictions of the \gp{} will be inaccurate. This can be seen in the figure for the KDE results (green diamonds) by their mean being away from zero and the large spread of their estimates (the outlier has a value of approximately $0.2$). However, due to the nature of this method, high values of $\widehat{S^H_{X^k}}$ are measured because the predicted output values are the values of the prediction mean function, which is a continuous function. Hence, these predictions are located on a curve. Therefore, the values of $\widehat{S^H_{X^k}}$ for the MST-based estimator are too large to be visible in this plot for the chosen values of $L$, except for $L=10^3$, where no emulated output is used. The reference result where we computed the FMST-based \si{} on the full data without emulator gave a reasonable result (red pentagon). 

Similar results appear for estimates based on emulation by stochastic collocation (SC), where we used collocation samples based on the uniform distribution. A function is fit through the data while no functional relation between input and output exists. Therefore, high values of $\widehat{S^H_{X^k}}$ are measured, which are not shown in the plot.

We summarize these results as follows: when an emulator (either \gp{} or SC) is used to augment the data for \sa{}, positive values of $\widehat{S^H_{X^k}}$ are found because the emulator is designed to fit a functional relation between input and output. The ``sample" method does not suffer from this problem. However, this is a very specific test case in which sample-based estimators are preferred over ones which use an augmented dataset.

\subsection{Analytic test case}\label{ssec:analytic}
\begin{figure}[ht!]
	\centering
	\includegraphics[width=0.5\textwidth]{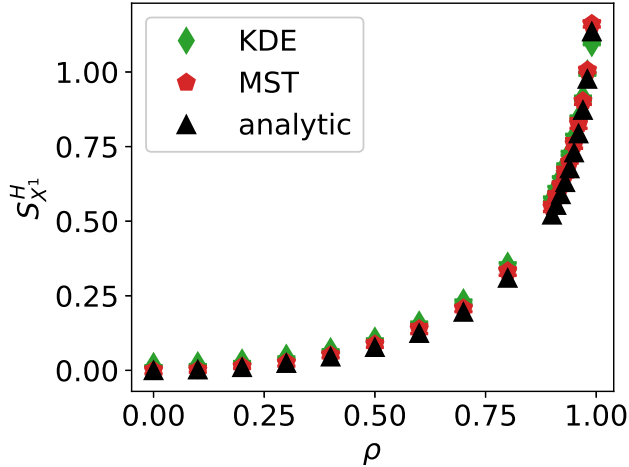}
	\caption{Computed values for the sensitivity indices, analytic test case.}\label{fig:exactsis}
\end{figure}

We consider a small test case in which we can compute the sensitivity index analytically. The idea behind this is to compare the KDE and the MST method in case no emulator is used. We have
\begin{equation*}
\begin{bmatrix}X \\ Y\end{bmatrix} \sim N\left(\begin{bmatrix}0 \\ 0\end{bmatrix},\begin{bmatrix}1 & \rho \\ \rho & 1\end{bmatrix}\right).
\end{equation*}
We took $N=10^4$ and repeated the experiment $n_r=10$ times. The results are in Figure \ref{fig:exactsis}.
Except for $\rho=0.98$, the MST method outperforms the KDE method. Furthermore, the MST method is i) not dependent on parameter choices such as kernel and kernel bandwidth and ii) faster to compute. One also needs to take into account that the rule of thumb to choose the kernel bandwidth we used here is based on the assumption that the data comes from a normal distribution and, therefore, this kernel bandwidth is optimal in this test case. When the underlying distribution is not normal, this heuristic may not be optimal. 

\newpage
\subsection{Ishigami function}\label{ssec:Ishi}
We now continue to a non-trivial synthetic test case, of which the test function is from Ishigami \& Homma \cite{3Ish90}. This output function is defined by
\begin{equation*}
	G(x,y,z|a,b) = (1+bz^4)\sin(x) + a\sin^2(y)
\end{equation*}
on the domain $[-\pi,\pi]^3$ (dimension $d=3$). We will use the well-known choice $a=7$, $b=0.1$ in accordance with \cite{2Cre09}.

Two types of input data are constructed for this test case. One is uniformly distributed and consists of $N=10^3$ samples on the domain of the output function. The other is the empirical copula of a multivariate normal distribution on the same domain, which is given by 
\begin{equation*}
Z = N\left(\begin{bmatrix}	0 \\ 0 \\ 0	\end{bmatrix},\begin{bmatrix} 1 & 0.8 & 0.5 \\ 0.8 & 1 & 0.8 \\ 0.5 & 0.8 & 1\end{bmatrix}\right),
\end{equation*}
such that
\begin{equation*}
X = -\pi + 2\pi\cdot F(Z),
\end{equation*}
with $F$ the cumulative distribution function of the marginal distributions (which is distributed as $N(0,1)$). 

For reference, we compute both the KDE-based and MST-based estimate on a larger dataset (with $N=10^5$) for comparison. Scott's rule \cite{4Sco92} is used for the kernel bandwidth.

In the numerical experiments, we first compute, depending on the dataset, a Latin hypercube sample (LHS) or Monte Carlo sample (MCS) of size $L=\{30,50,100,200\}$ and combine it with KDE. For this data, we computed (\ref{eq:H}). Then, we fit a \gp{} with Gaussian kernel to these samples, where the length scales have been estimated by maximum likelihood estimation. Now, we can proceed with KDE on $(X_+,Y_+)$, in which we include the choice $Y_{L^+}=\mu(X_{L^+})$ (Equation \ref{eq:est1}). We obtain one estimate for $\overline{H}_{X^k,f}^{(L+L_+)}$ for each repetition of the experiment and thereby one value of $|\overline{H}_{X^k,f}^{(L+L_+)}-\overline{H}_{X^k,f}^{(N)}|\approx |\hat{S}_{X^k}-S_{X^k}|$ which is used as measure of convergence. 
\added{In a similar way, we can proceed with the MST method on $(X^+,Y^+)$ with $Y_{L^+}=\mu(X_{L^+})$. Finally, the SC method, based on the uniform distribution and combined with KDE, is used for comparison. Note we showed earlier that KDE has a larger bias than MST, but we look mainly at the convergence here. }

\begin{figure}[ht!]
	\centering
	\includegraphics[width=0.5\textwidth]{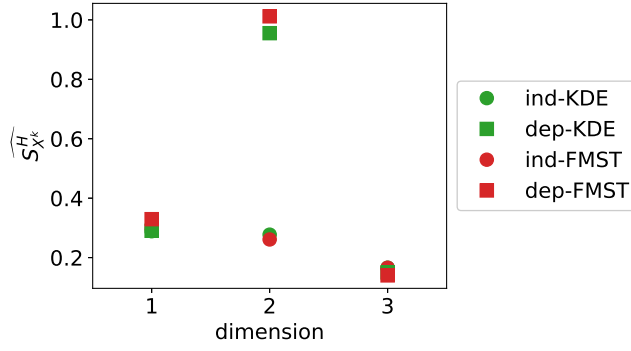}
	\caption{Computed values for the sensitivity indices per variable, Ishigami function.}\label{fig:sis}
\end{figure}

The computed reference values of the \sis{} are shown in Figure \ref{fig:sis}. 
The values for independent and dependent data are close to each other for variables $1$ and $3$, while they are far apart for variable $2$, which is due to the dependency. 

We will first show the results for the independent data, followed by the results for the dependent data. We start with determining the goodness-of-fit of the \gp{} by performing $k$-fold cross-validation (CV) with $k=10$ and compute the coefficient of determination
\begin{equation*}
R^2=1-\frac{SS_{res}}{SS_{tot}} = 1-\frac{\frac{1}{L}\sum_l (\hat{Y}_l-Y_l)^2}{\frac{1}{L}\sum_{l} (Y_l-\bar{Y})^2},
\end{equation*}
where $\hat{Y}_l$ are the CV predictions for $Y_l$ and $\bar{Y}=\frac{1}{L}\sum_{l}Y_l$. 
\begin{figure}[ht!]
	\centering
	\includegraphics[width=0.3\textwidth]{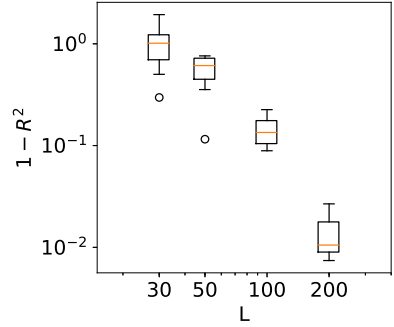}
	\caption{Cross-validation results showing the quality of the \gp{}, independent input data.}\label{fig:cv}
\end{figure}
The values for $R^2$ for independent data are given in Figure \ref{fig:cv} and we see its values are near zero for higher values of $L$. For $L=30$ and $L=50$, this fraction can become larger than $1$. In this case, the fit is worse than a constant function. Note that here, the \gp{} is not fit well, while this is the case for the higher values of $L$.

\begin{figure}[ht!]
	\centering
	\begin{subfigure}{0.3\textwidth}
		\includegraphics[height=0.17\textheight,trim=0mm 0mm 45mm 0mm,clip]{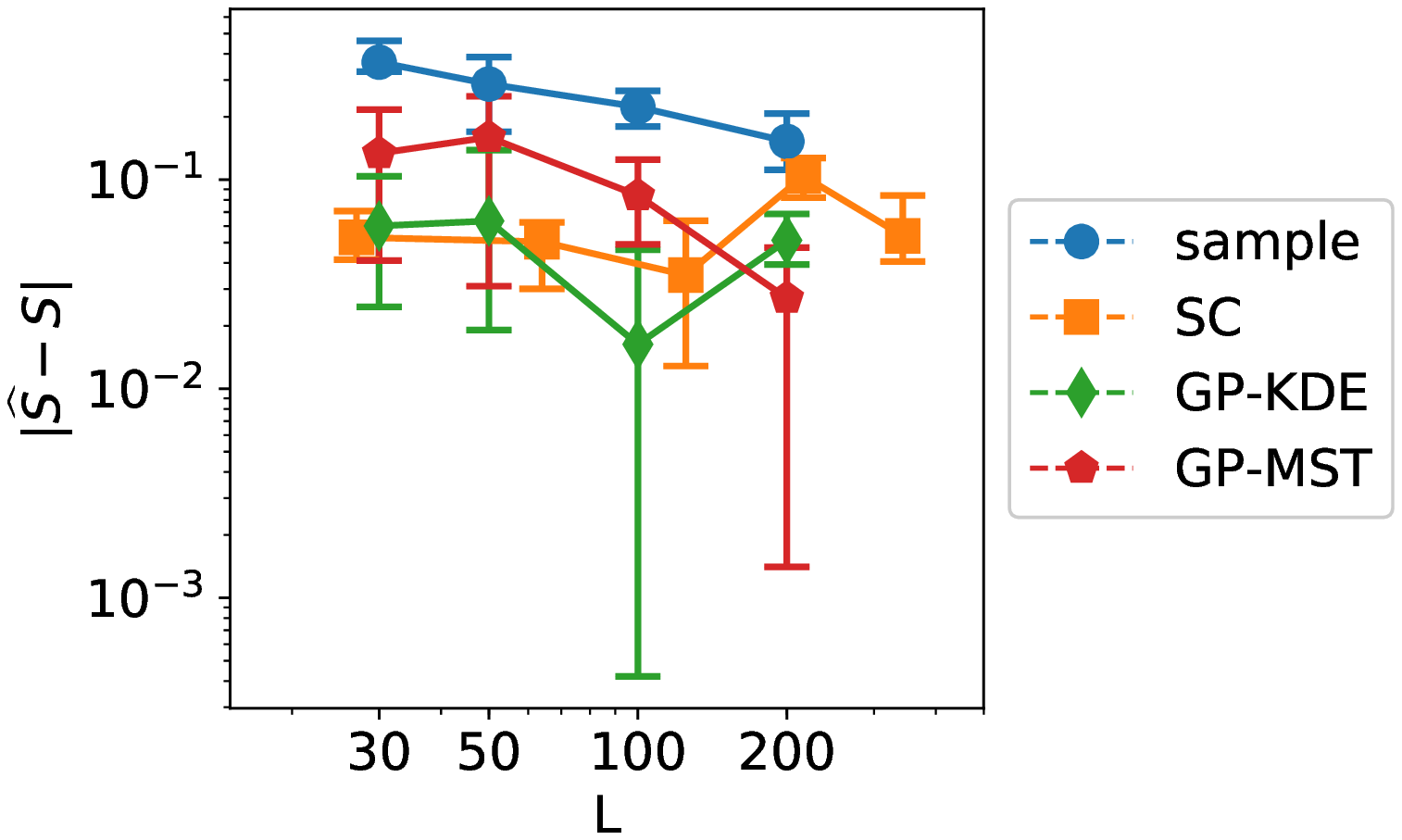}
	\end{subfigure}
	\begin{subfigure}{0.3\textwidth}
		\includegraphics[height=0.17\textheight,trim=0mm 0mm 45mm 0mm,clip]{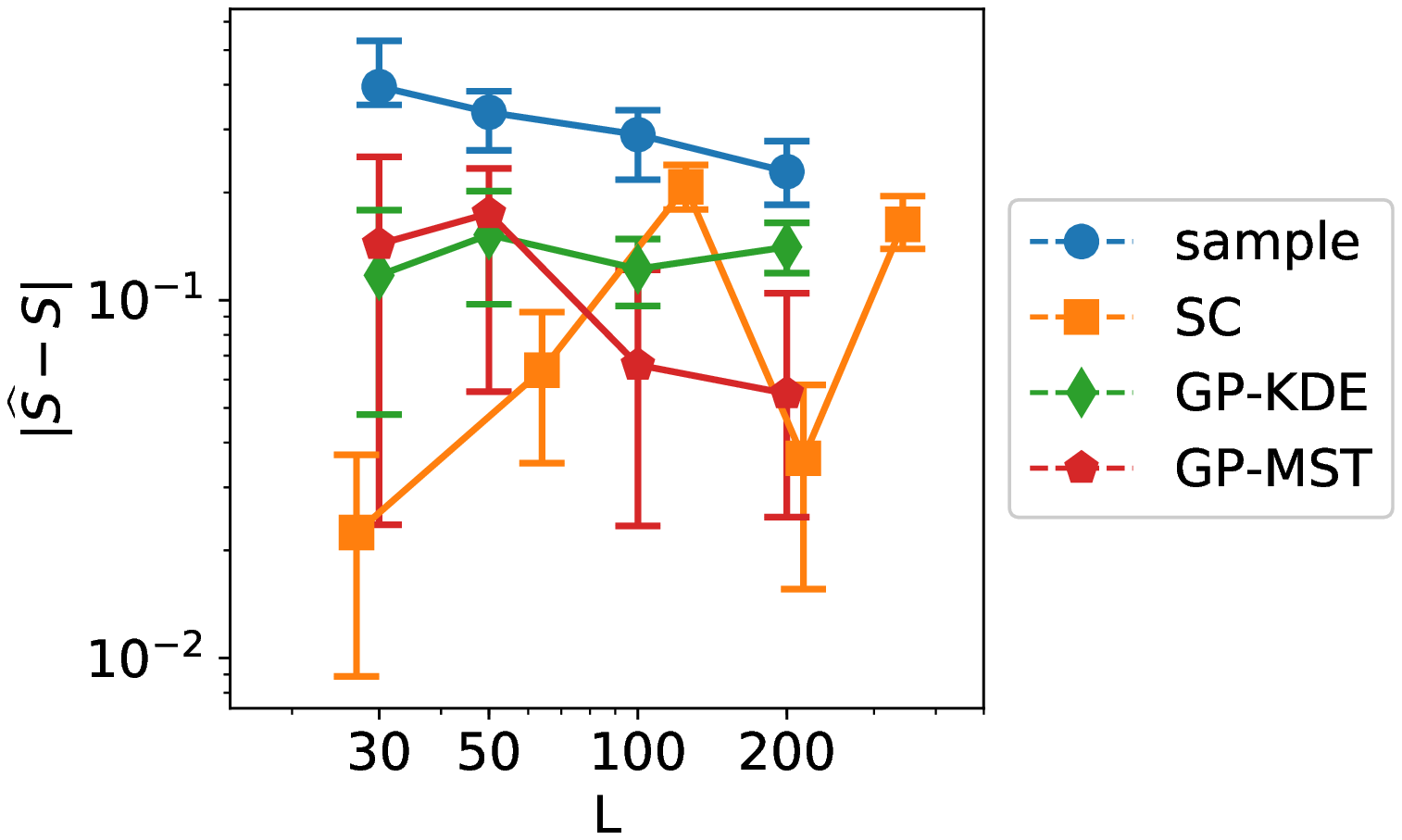}
	\end{subfigure}
	\begin{subfigure}{0.38\textwidth}
		\includegraphics[height=0.17\textheight]{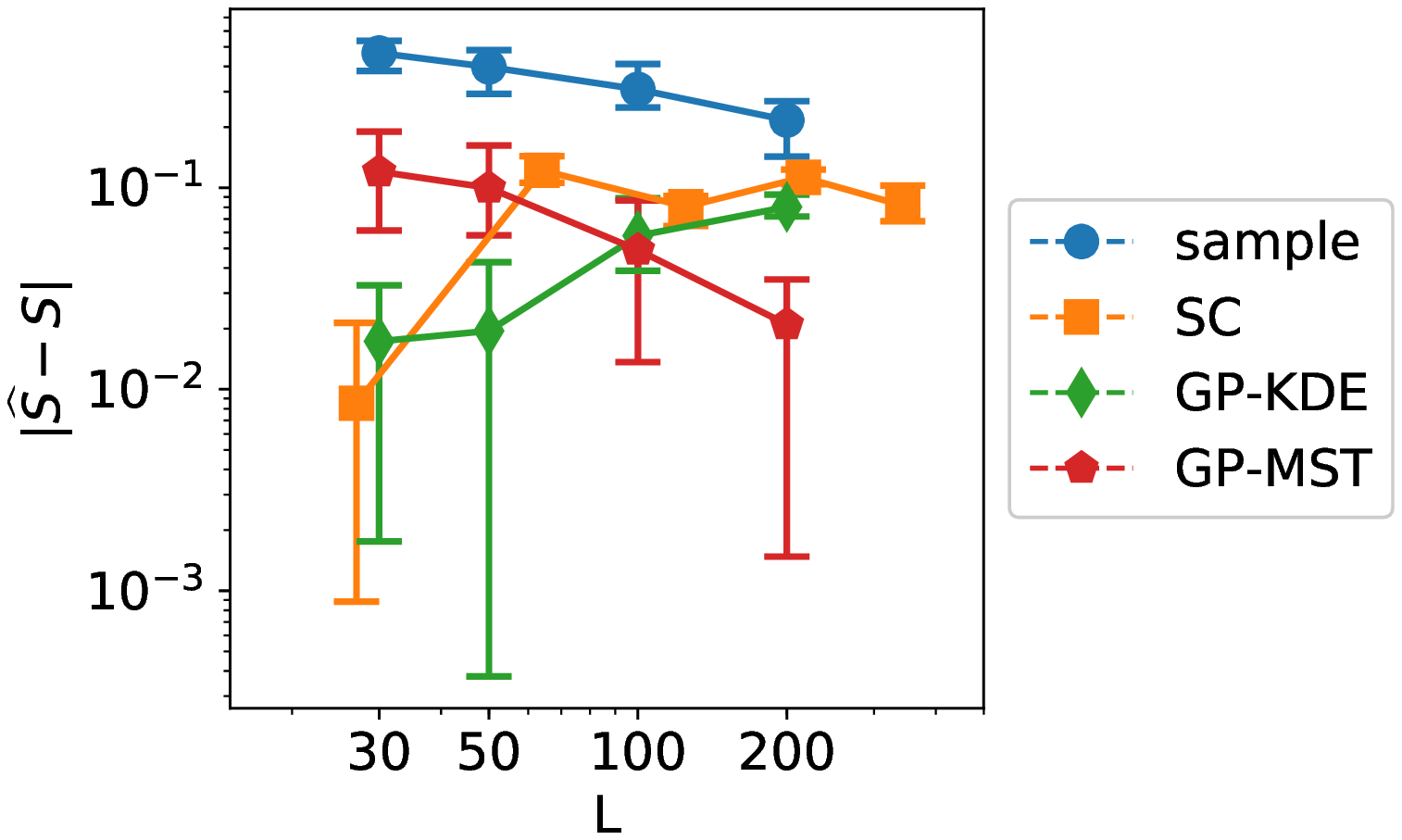}
	\end{subfigure}
	\caption{Convergence of the estimates for Hellinger divergence, independent input data. }\label{fig:convKL}
\end{figure}

Figure \ref{fig:convKL} shows the convergence of the estimates, where ``sample'' indicates the KDE is based on only $L$ samples, ``SC" indicates stochastic collocation is used (combined with KDE), ``GP-KDE" is based on (\ref{eq:est1}) and ``GP-MST" is based on (\ref{eq:SImst}). From left to right, variables 1 to 3 are shown. This will also be the case for all similar figures in this section. 
In this figure, we see several trends. First of all, the samples perform worse than the methods which use augmentation of the dataset. Second, we see the results for SC are not robust and their errors do not decrease in general for increasing $L$. Third, we see that GP-MST shows in general decreasing errors for increasing $L$. 

\begin{figure}[ht!]
	\centering
	\includegraphics[width=0.3\textwidth]{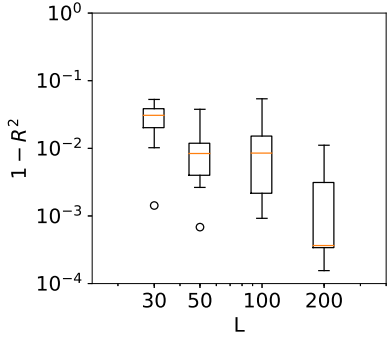}
	\caption{Cross-validation results showing the quality of the \gp{}, dependent input data.}\label{fig:cvdep}
\end{figure}

\begin{figure}[ht!]
	\centering
	\begin{subfigure}{0.3\textwidth}
		\includegraphics[height=0.17\textheight,trim=0mm 0mm 45mm 0mm,clip]{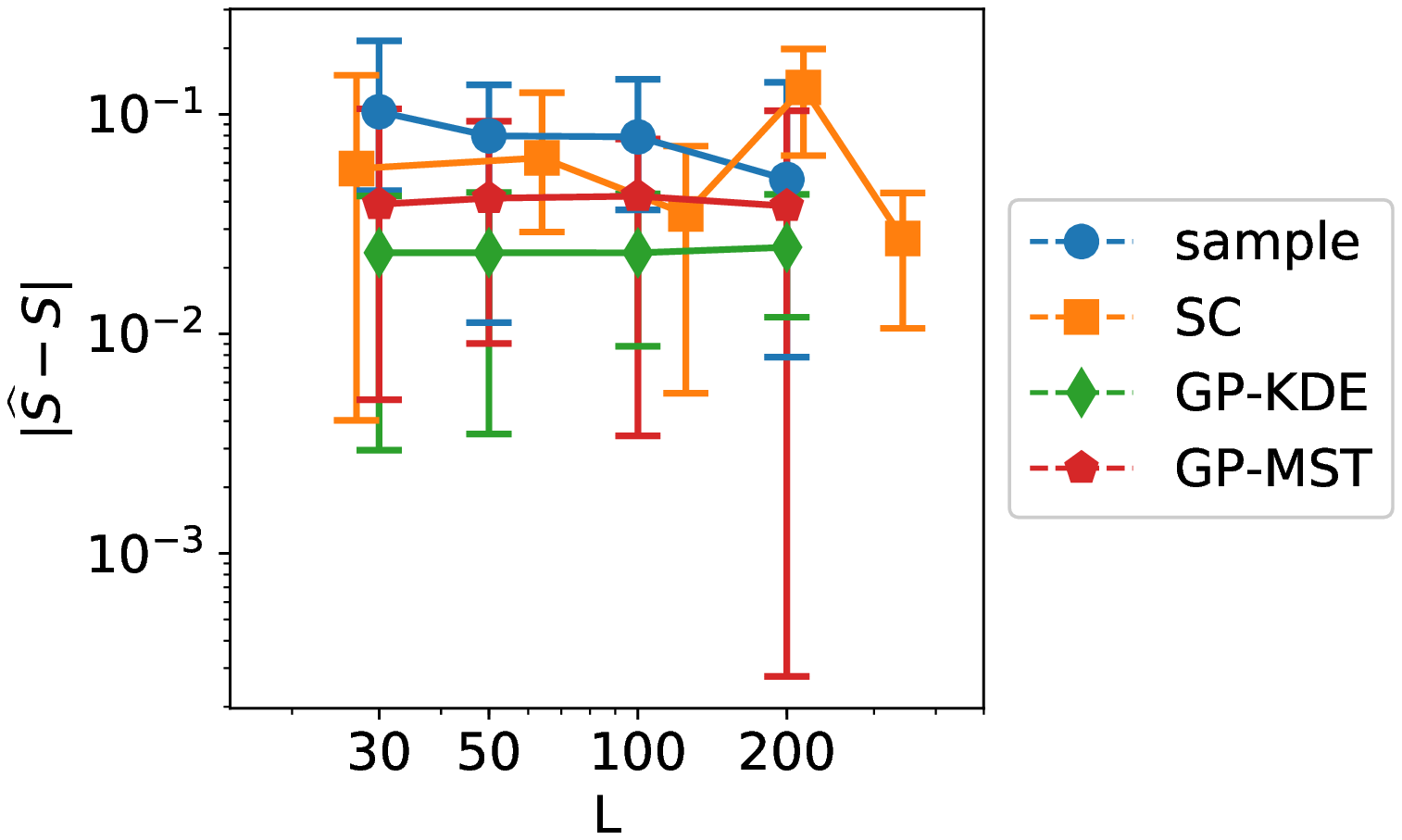}
	\end{subfigure}
	\begin{subfigure}{0.3\textwidth}
		\includegraphics[height=0.17\textheight,trim=0mm 0mm 45mm 0mm,clip]{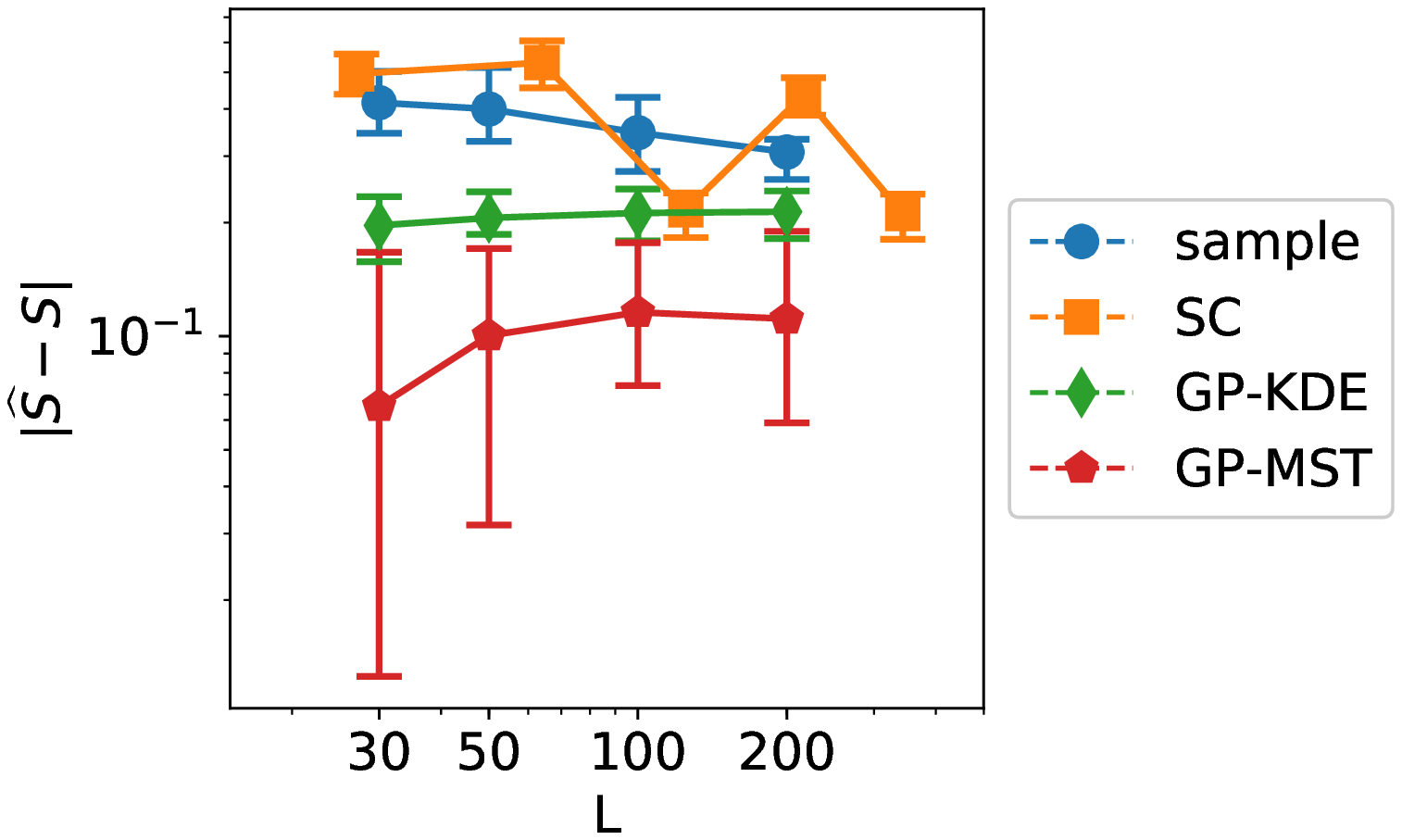}
	\end{subfigure}
	\begin{subfigure}{0.38\textwidth}
		\includegraphics[height=0.17\textheight]{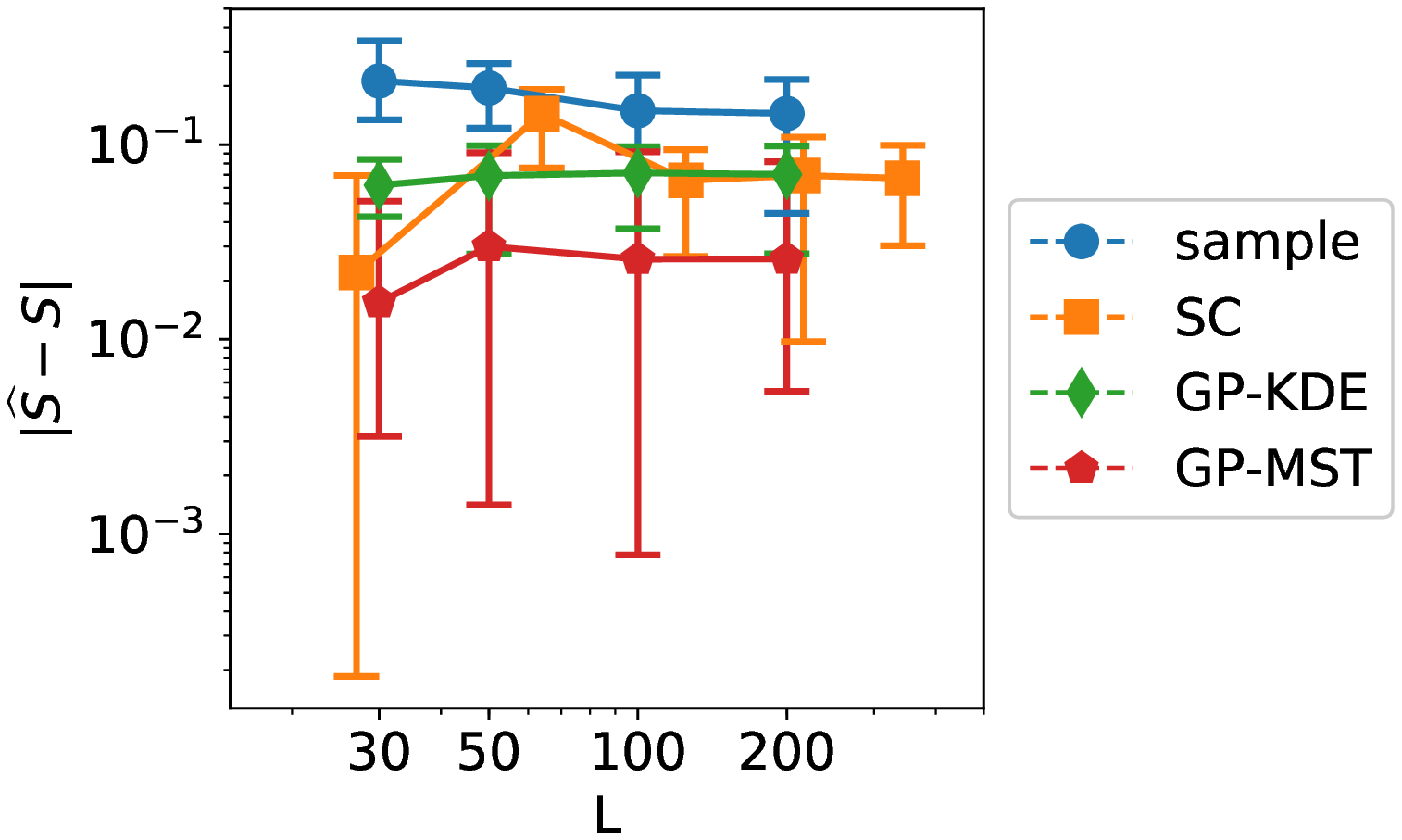}
	\end{subfigure}
	\caption{Convergence of the estimates for Hellinger divergence, dependent input data.}\label{fig:convKLdep}
\end{figure}

The results for dependent data are shown in Figures \ref{fig:cvdep} and \ref{fig:convKLdep}. Note that LHS is not an appropriate sampling method because the data is dependent, therefore, Monte Carlo sampling is used instead. Furthermore, SC is here also not completely suitable because the input distribution is dependent. The results are similar to previous experiments, although the cross-validation results imply the \gp{} for this data has been fit better. Another observation is that GP-MST outperforms the other methods for variables 2 and 3, while it is not really worse than GP-KDE for variable 1. Overall, the \gp{}-based methods outperform the other methods.

\subsection{Piston function}\label{ssec:Piston}
\begin{figure}[ht!]
	\centering
	\includegraphics[width=0.45\textwidth]{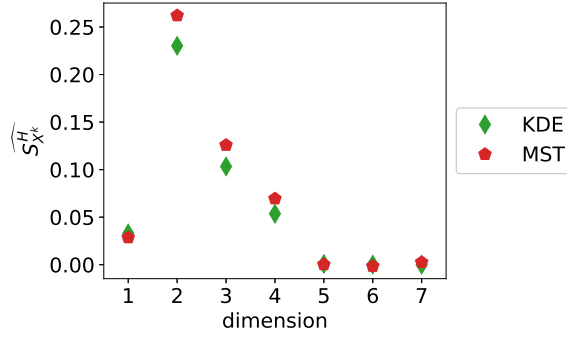}
	\caption{Computed values for the sensitivity indices per variable, Piston function.}\label{fig:sisP}
\end{figure}

\begin{figure}[ht!]
	\centering
	\includegraphics[width=0.3\textwidth]{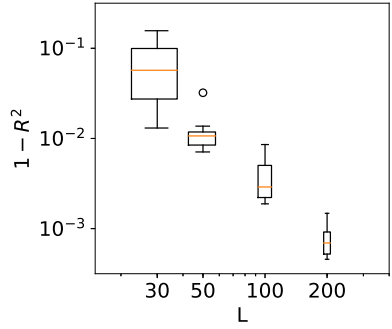}
	\caption{Cross-validation results showing the quality of the \gp{}, Piston function.}\label{fig:cvP}
\end{figure}

\begin{figure}[ht!]
	\centering
	\begin{subfigure}{0.24\textwidth}
		\includegraphics[height=0.14\textheight,trim=0mm 0mm 45mm 0mm,clip]{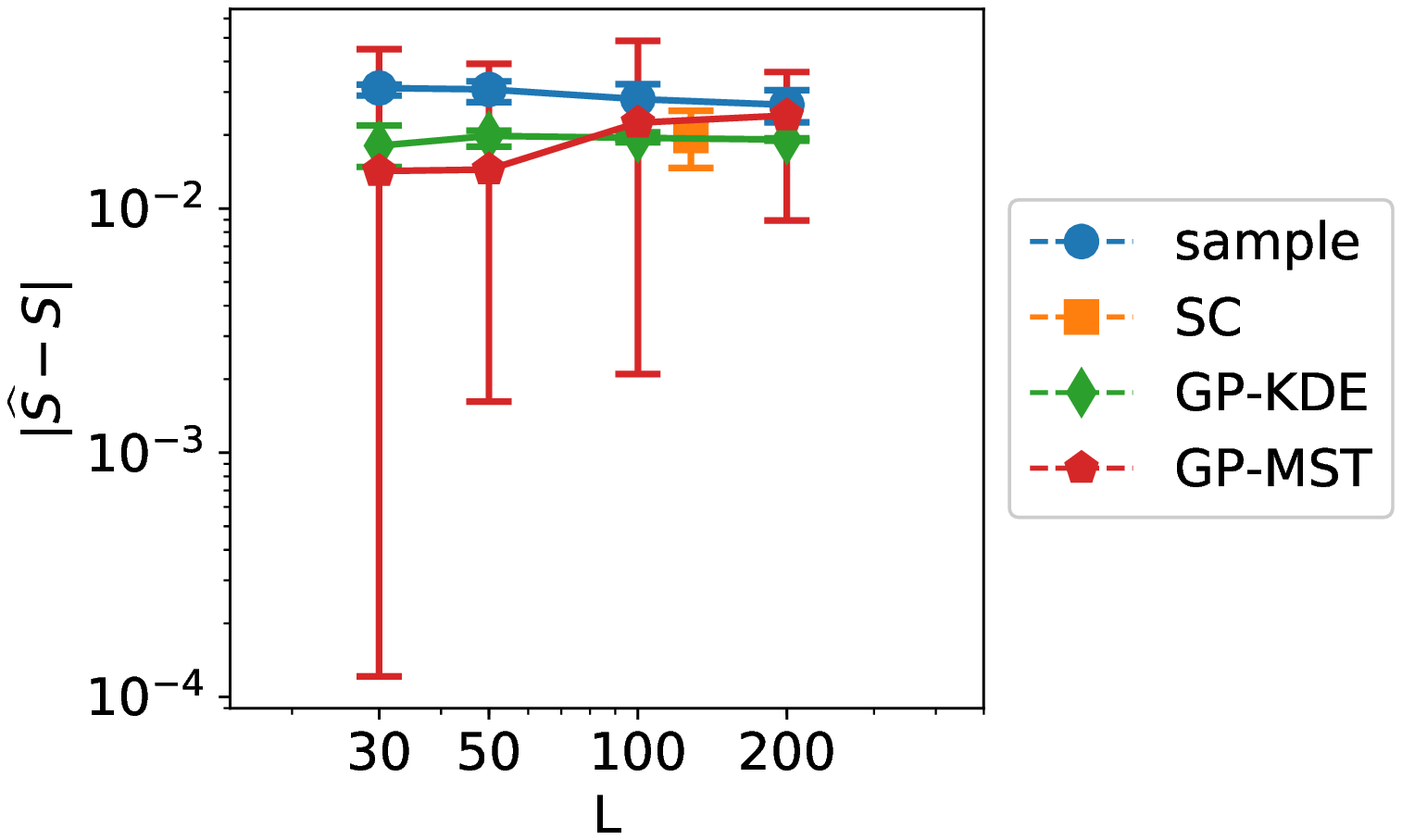}
	\end{subfigure}%
	\begin{subfigure}{0.26\textwidth}
		\includegraphics[height=0.14\textheight,trim=0mm 0mm 43mm 0mm,clip]{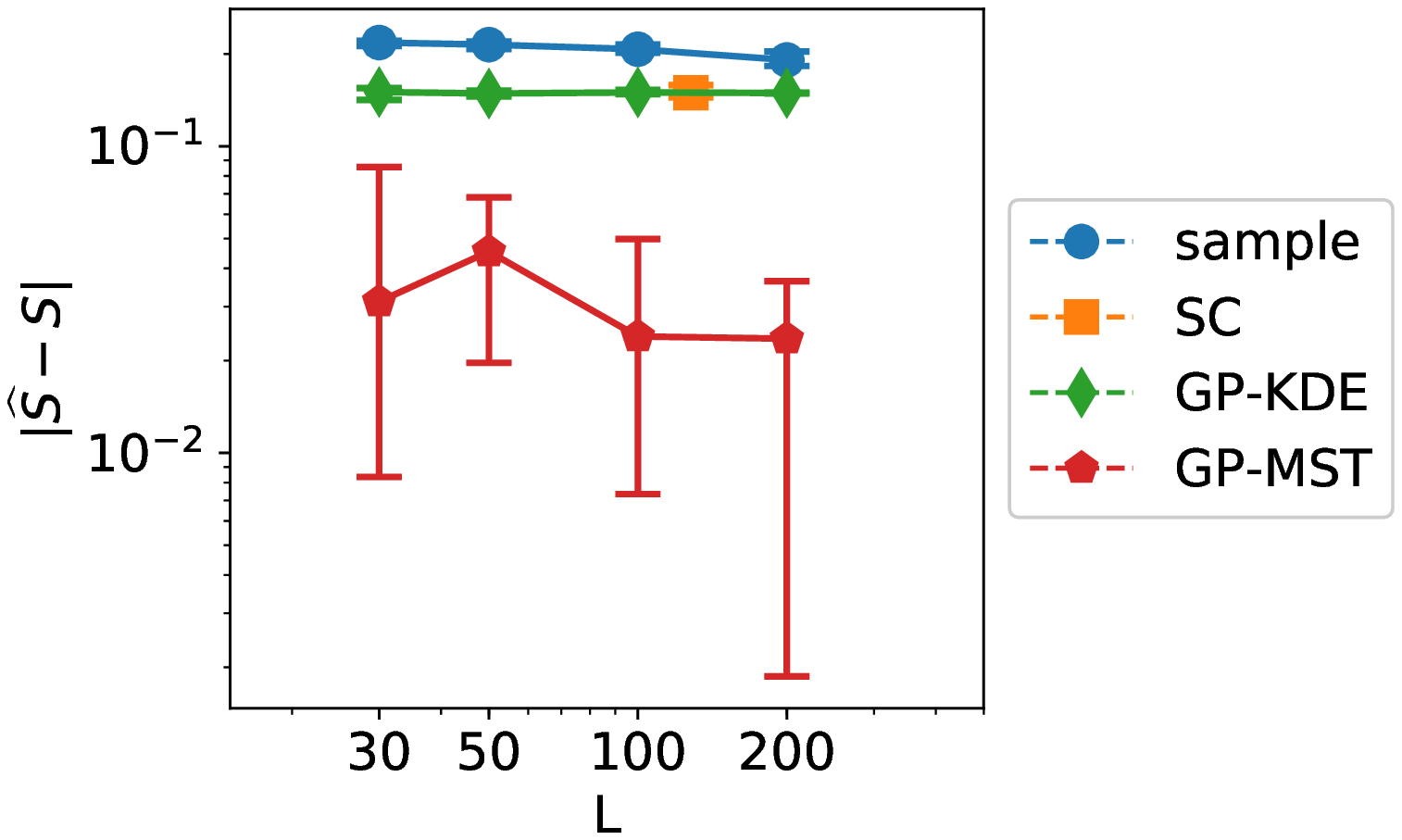}
	\end{subfigure}%
	\begin{subfigure}{0.26\textwidth}
		\includegraphics[height=0.14\textheight,trim=0mm 0mm 43mm 0mm,clip]{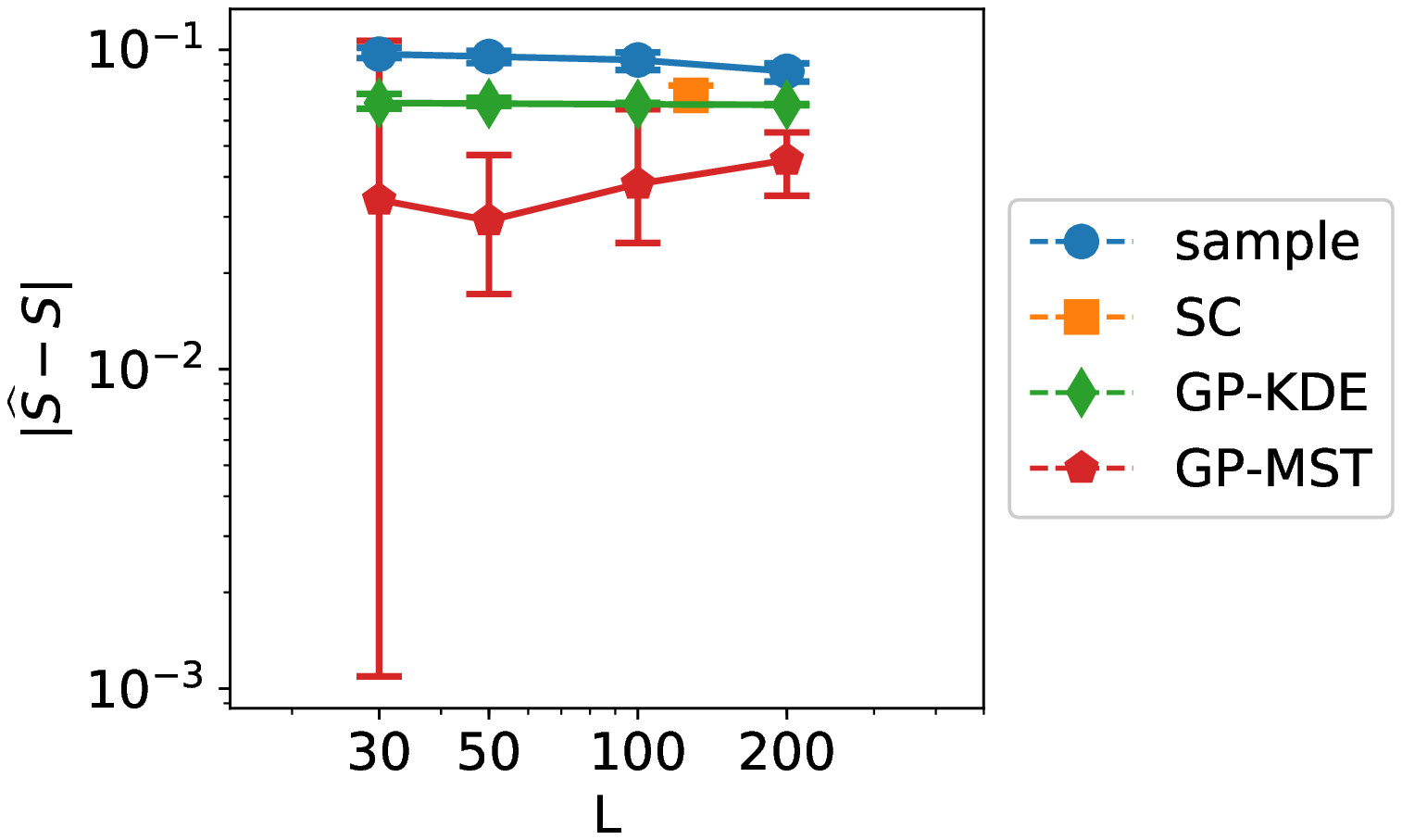}
	\end{subfigure}%
	\begin{subfigure}{0.24\textwidth}
		\includegraphics[height=0.14\textheight,trim=0mm 0mm 45mm 0mm,clip]{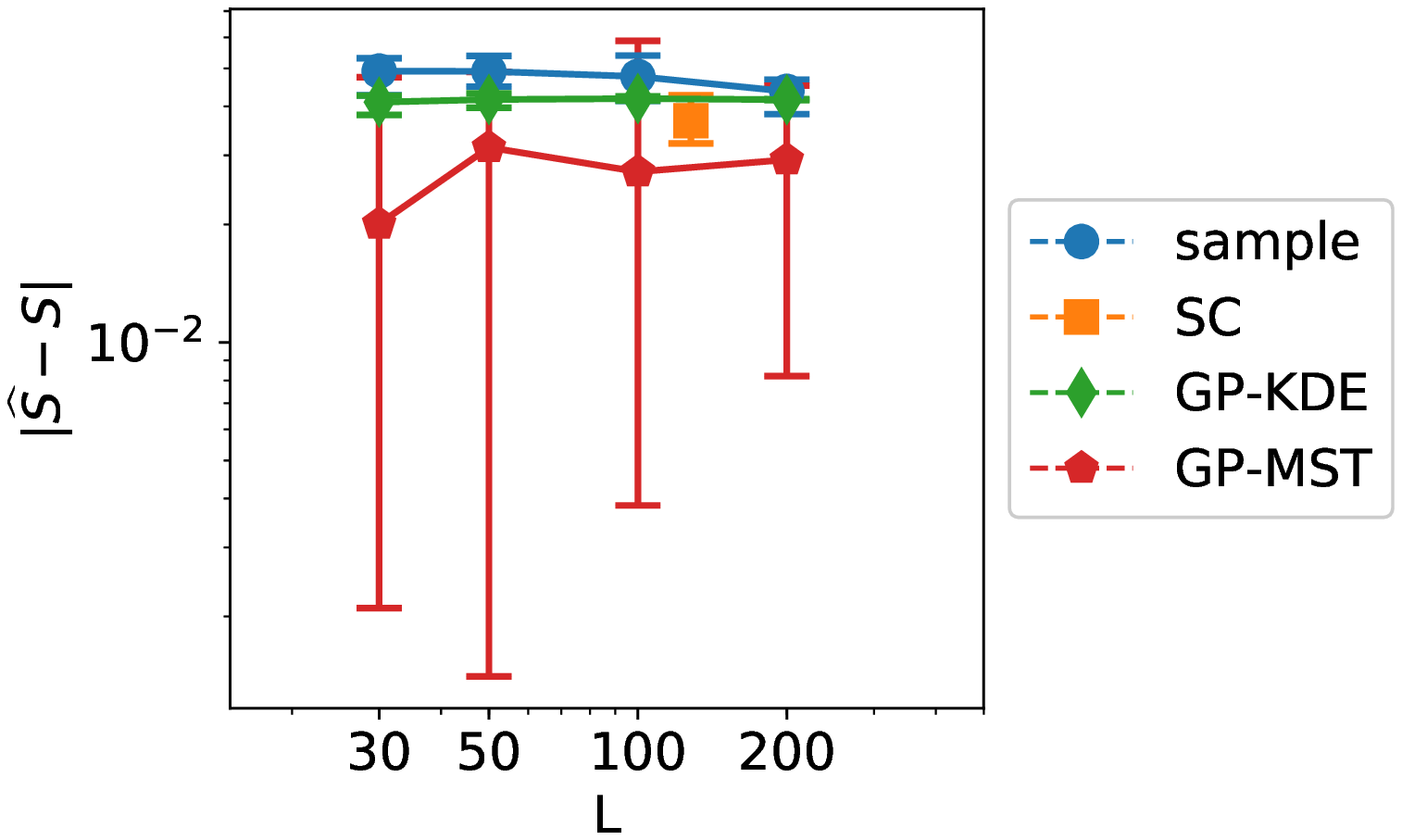}
	\end{subfigure}\\
	\begin{subfigure}{0.24\textwidth}
		\includegraphics[height=0.14\textheight,trim=0mm 0mm 45mm 0mm,clip]{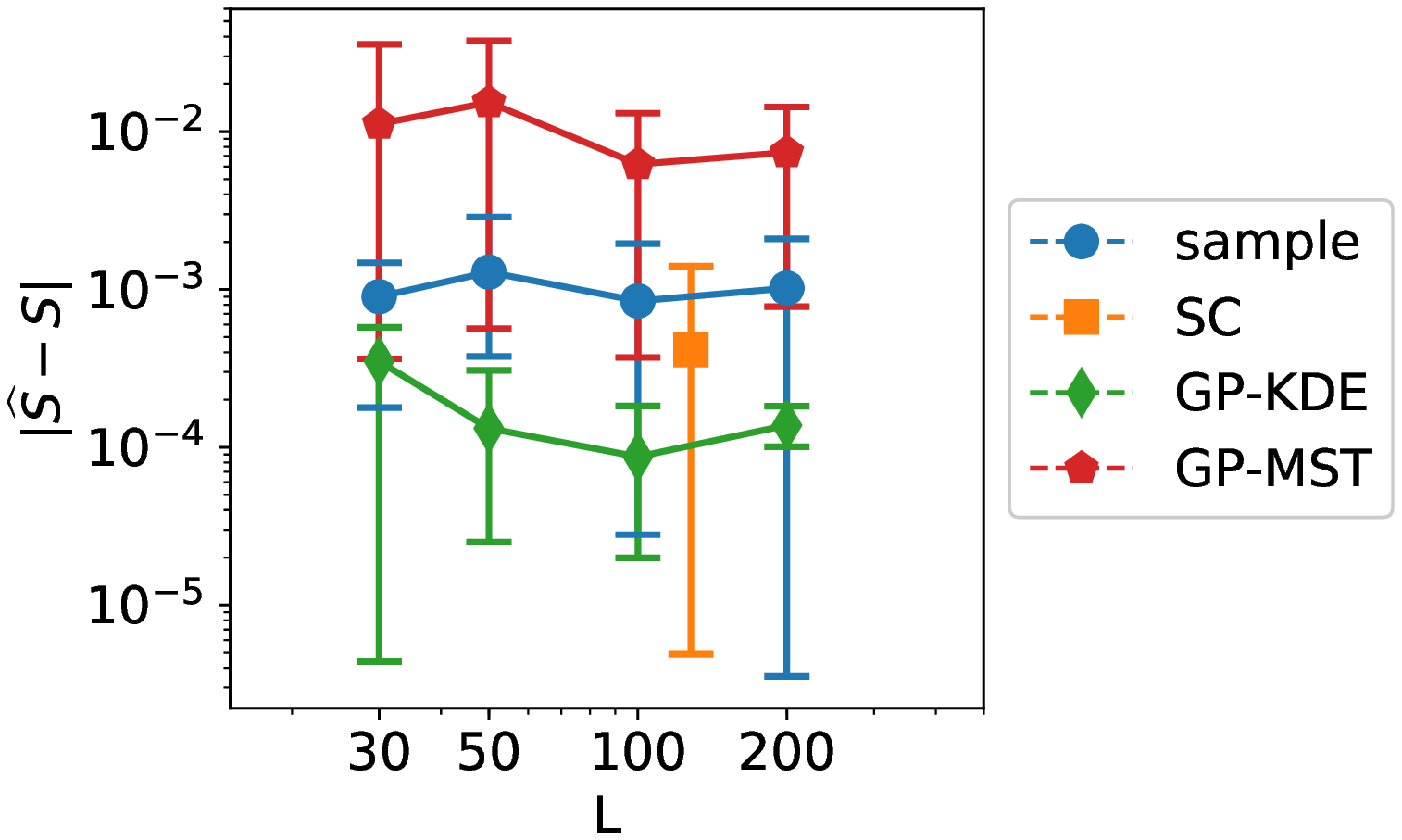}
	\end{subfigure}%
	\begin{subfigure}{0.24\textwidth}
		\includegraphics[height=0.14\textheight,trim=0mm 0mm 45mm 0mm,clip]{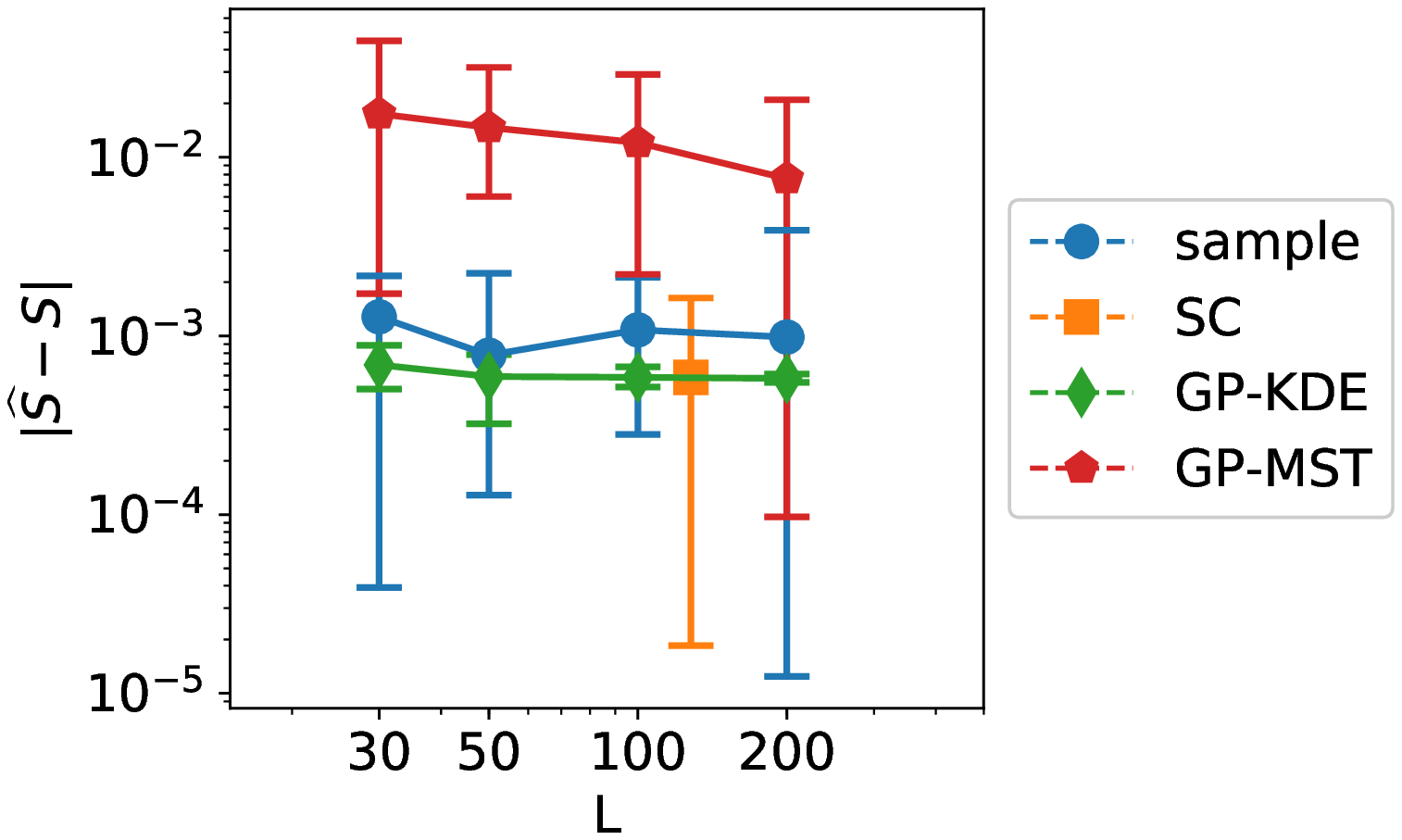}
	\end{subfigure}%
	\begin{subfigure}{0.3\textwidth}
		\includegraphics[height=0.14\textheight]{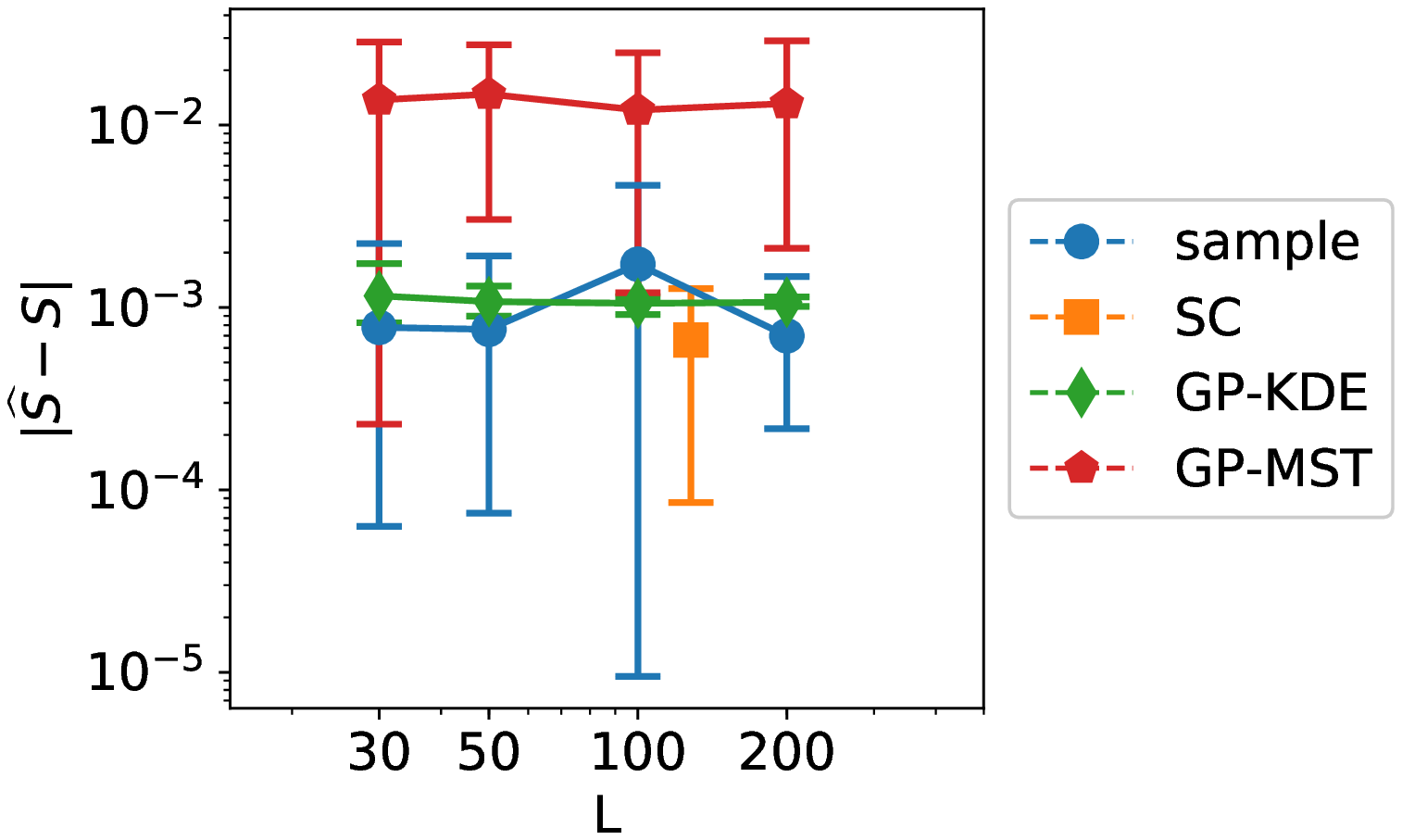}
	\end{subfigure}%
	\caption{Convergence of the estimates for Hellinger divergence, Piston function.}\label{fig:convKLP}
\end{figure}

We also tested a higher-dimensional test case with independent uniformly distributed input variables. In this case, the output function is defined by the Piston function from \cite{4Ken98}. The output here is the cycle time of a piston, as given by
\begin{align*}
C(\mathbf{x}) = & 2\pi\sqrt{\frac{M}{k+S^2\frac{PV}{T_0}\frac{T_a}{V^2}}}, \\
V & = \frac{S}{2k}\left(\sqrt{A^2+4k\frac{PV}{T_0}T_a}-A\right) \\
A & = PS + 19.62M-\frac{kV}{S} \\
\mathbf{x} & = (M,S,V,k,P,T_a,T_0), 
\end{align*}
of which the input ranges are given in Table \ref{tab:vars}.
\begin{table}
	\centering
	\caption{Input variables for the Piston function.}\label{tab:vars}
	\begin{tabular}{r|l}
		\hline
		Symbol and range & Explanation \\
		\hline
		$M \in [30,60]$ & piston weight (kg) \\
		$S \in [0.005,0.020]$ & piston surface area (m$^2$)\\
		$V \in [0.002,0.010]$ & initial gas volume (m$^3$) \\
		$k \in [1000,5000]$ & spring coefficient (N/m) \\
		$P \in [90000,110000]$ & atmospheric pressure (N/m$^2$) \\
		$T_a \in [290,296]$ & ambient temperature (K) \\
		$T_0 \in [340,360]$ & filling gas temperature (K) \\
	\end{tabular}
\end{table}
For numerical reasons, the data of size $N=10^3$ is generated and processed on the unit hypercube: it is only transformed to the input ranges to obtain the output values. The \sis{} as computed on a larger dataset of size $N=10^5$ are given in Figure \ref{fig:sisP}. The values for KDE and MST differ, although Section \ref{ssec:analytic} indicates the MST results are more accurate.
The cross-validation results are in Figure \ref{fig:cvP}. These last results show the \gp{} has been fit well for $L\geq 50$ and therefore we can continue with the remaining results.
The results for the convergence are in Figure \ref{fig:convKLP}. From left to right, top to bottom, variables 1 to 7 are shown. 

We see clear differences between variables 1-4 on one hand and variables 5-7 on the other hand. This is due to variables 5-7 for which $\widehat{S^H_{X^k}}$ is near zero. As indicated in Section \ref{ssec:random}, methods which make use of a fit perform badly in this case. For variables 1-4, GP-MST clearly outperforms the others. The SC method has only been performed with $2$ collocation points for each dimension, which led to $L=2^7=128$. Increasing to $3$ would give us $3^7=2187$ collocation points, which is higher than the number data points in the used dataset.

\subsection{Recommendation}\label{ssec:rec}
The \gp{}-based methods in general outperform the sample-based method and stochastic collocation, except when the value of the \si{} is (near-)zero. However, usually one is interested in ordering the input variables based on the \sis{} rather than obtaining their values very precisely. Input variables with values of the \si{} near zero are usually considered unimportant and in that case, it is also not very important to estimate the value of zero very precisely. We therefore advise to use the GP-MST method, wherein the available samples $(X_L,Y_L)$ are augmented to $(X_+,Y_+)$ by a \gp{}, on which the \si{} (\ref{eq:sif}) is computed for the Hellinger distance by the minimum spanning tree method.

\section{Direct \sis{}}\label{sec:dsis}
We note that the \sis{} as described by \cite{4Dav14} are total \sis{}, which include both direct and indirect effects. Direct effects measure the effect of one input variable only, while indirect effects contain the effect of the other variables due to possible dependencies in the input variables. The indirect effect is the difference between the total and the direct effect. To illustrate this, consider an example where $X=(X^1,X^2)$ follows a bivariate normal distribution with means $0$, variances $1$ and covariance $\rho>0$ (so that $X^1$ and $X^2$ are dependent), while $u(x_1,x_2)=2x_1$. Then the direct effect of $X^2$ is zero, while its total effect is positive (because $u(X^1,X^2)$ and $X^2$ are dependent through $X^1$). Hence, the indirect effect of $X^2$ is positive as well. Our goal is now to find a measure for the direct effects, i.e., without the effects of the mutual input dependencies. 

Although useful, total \sis{} do not tell the complete story. While an input variable may be completely irrelevant for the value of the output, it may have a positive sensitivity index due to a dependency with a relevant input. The relevant input variable would then be called a confounder. An example of this is wave height for a computational model of offshore wind energy: although the waves have nothing to do with the power output, they are linked to each other via the wind speed with which they have a dependency. To get rid of this effect, we need to construct indices which measure the effect of only one input variable, without effects due to dependencies in the input. It is in this case necessary to remove the dependencies from the input.

For variance-based sensitivity indices, a distinction is made between first-order, higher-order and total sensitivity indices \cite{3Sob01}. In first-order indices, one only measures the effect of varying one variable alone, where in higher-order indices multiple variables are varied at the same time. Because the number of second-order sensitivity indices grows as $d(d-1)/2$ with $d$ the number of input variables, and the total number of sensitivity indices is $2^d-1$, usually not all of them are computed. Instead, one computes the first-order and total sensitivity indices. 

In a similar fashion to first-order indices, we define direct \sis{}, which measure the effect of varying one variable only. The direct indices then measure the direct effects, while the total indices measure the combination of direct and indirect effects, which also includes effects due to dependencies in the input.

\subsection{Theory}
The starting point of these new indices is the same divergence-based index as before, namely (\ref{eq:sih}). We repeat (\ref{eq:sih}) here for convenience,
\begin{equation*}
S^H_{X^k} = 2 - 2 \iint_{\mathbb{R}^2}\sqrt{p_{X^k}(x)p_Y(y)p_{X^k,Y}(x,y)}\mathrm{d}y\mathrm{d}x.
\end{equation*}
Now, note that
\begin{equation*}
Y = u(X^1,\ldots,X^d) = u(X),
\end{equation*}
with $u(\cdot)$ being the model used to obtain the output $Y$, hence, both $p_Y(y)$ and $p_{X^k,Y}(x,y)$ depend in theory on all input variables $X^k$. Hence, if we remove the dependencies between the input variables, then these probability distributions change as well. This removal is done by applying a permutation operator $\Pi$, which is defined on a dataset $X$ in such a way that
\begin{equation*}
\Pi(X) = (\Pi(X^1),\ldots,\Pi(X^d)),
\end{equation*}
with
\begin{equation*}
CDF(\Pi(X^k)) = CDF(X^k), \quad \Pi(X^i) \perp \Pi(X^j) \, \text{ for all } i\text{,}j, \, \, \, \, i\neq j.
\end{equation*}
Hence, this operator keeps the marginal distributions the same, but it removes all dependencies ($\perp$ here denotes statistical independence). The implementation of this operator is detailed at the end of this section.

Now, we create a permuted version of our dataset $X$, being $\Pi(X)$. For this dataset, we can define the direct \si{} by
\begin{equation*}
S^H_{D,X^k} = 2 - 2 \iint_{\mathbb{R}^2}\sqrt{p_{\Pi(X^k)}(x)p_{Y}(y)p_{\Pi (X^k),Y}(x,y)}\mathrm{d}y\mathrm{d}x.
\end{equation*}
The output $Y=u(X)$ can be replaced by
\begin{equation*}
\widetilde{Y}=\tilde{u}(\Pi(X)),
\end{equation*}
in which $\tilde{u}(\cdot)$ denotes the \gp{} $\widetilde{G}_{\{X_L,Y_L\}}(\cdot)$ constructed earlier. This leads to the estimator
\begin{equation*}
\overline{S}^H_{D,X^k} = 2 - 2 \iint_{\mathbb{R}^2}\sqrt{p_{\Pi(X^k)}(x)p_{\widetilde{Y}}(y)p_{\Pi (X^k),\tilde{Y}}(x,y)}\mathrm{d}y\mathrm{d}x.
\end{equation*}

The problem now is how to define $\Pi(X)$. A naive implementation could be one in which for each variable, a random permutation of the values is performed. This is fast, but does not guarantee independence of the input variables after transformation. Also, the indexing of the permutations leads to a Latin hypercube design (LHD): each value from $1$ to $N$ (for $N$ data points in the dataset) is used only once. However, this does not guarantee all dependencies are removed. In Latin hypercube sampling (LHS), a comparable problem exists as equally probable subspaces can end up with a different number of sampling points. This is solved by orthogonal sampling \cite{4Owe92} or by using a maximin criterion \cite{3Joh90}. 

Inspired by this, we would like to generate an LHD of size $N$ in $d$ dimensions with the maximin criterion which puts the samples at the middle of each interval. This LHD is easily transformed to an indexing, which can be applied to the original data $X$ to obtain $\Pi(X)$. However, obtaining such an LHD is computationally very expensive because it contains an optimization step and is therefore not feasible for the problem sizes we are looking at.

An alternative to Latin hypercube sampling is quasi-Monte Carlo sampling, which generates data points from low-discrepancy sequences such as Halton's \cite[Chapter 3]{4Nie92} and Sobol's \cite{4Joe08}. In this way, we achieve the goal that the proportion of data points in a sequence falling into a subspace is nearly proportional to the probability measure of this subspace (the difference between them is the discrepancy). Hence, we achieve an approximately uniform distribution of data points over the unit hypercube, which means the dimensions are independent of each other. Furthermore, all values generated for a variable are unique, which means they can easily be transformed to the discrete hypercube $\{1,\ldots,N\}^d$. The transformed values can be used as an indexing for $X$ to obtain $\Pi(X)$. Because the data points generated by the sequence are uniform over the unit hypercube, they lead to an independent dataset when their transformed values are used as indexing. We use the Sobol sequences as described by \cite{4Joe08}.

\subsection{Ishigami}\label{ssec:dsisIshi}
We compute these direct \sis{} for the Ishigami test case of Section \ref{ssec:Ishi}. We split the results out to the KDE and the MST estimates. For each of them, we show the estimates of both the independent and dependent direct \sis{} and the spread therein for increasing $L$. We also compare them to the values of the total indices.

\begin{figure}[ht!]
	\centering
	\begin{subfigure}{0.48\textwidth}
		\includegraphics[height=0.195\textheight]{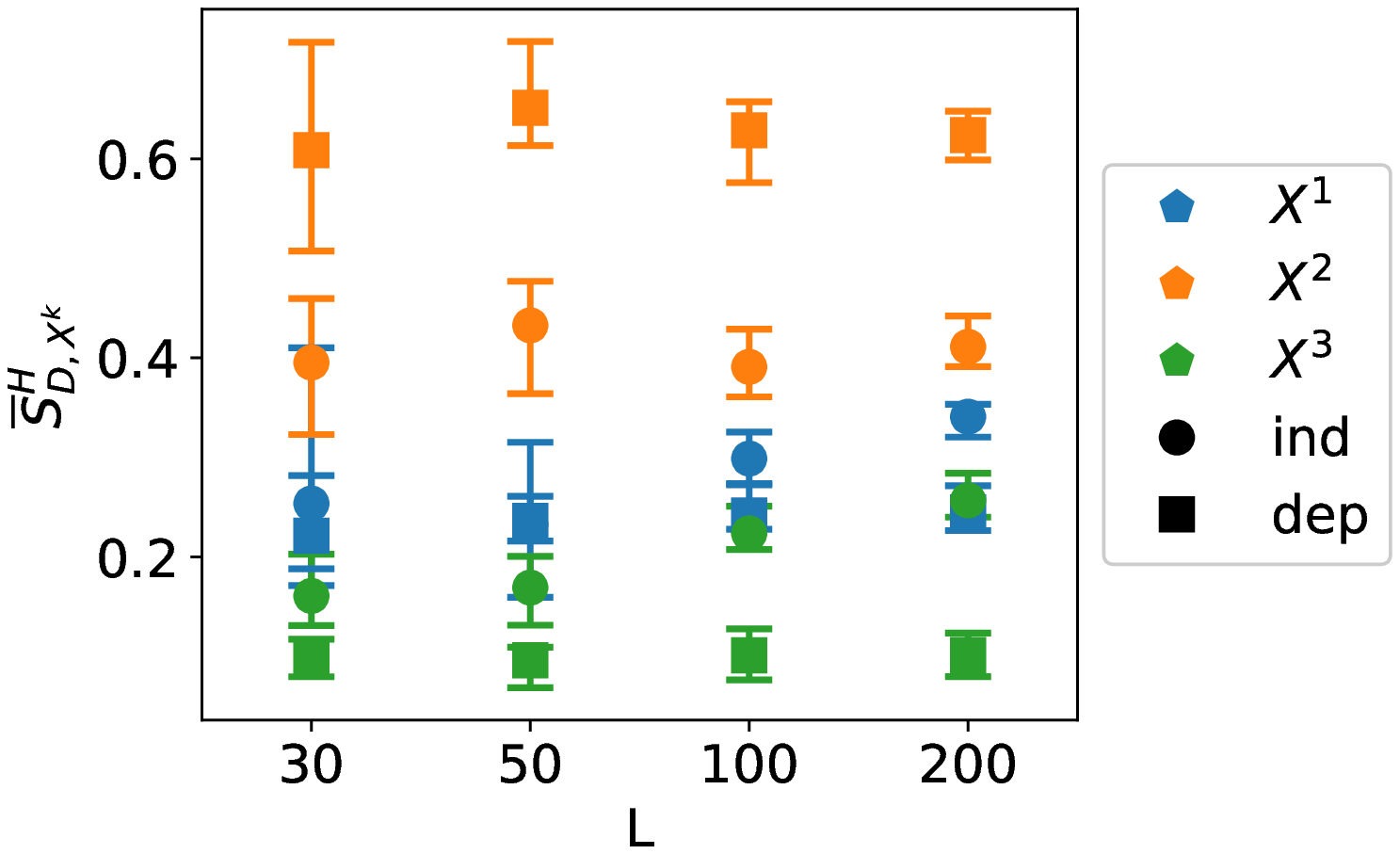}
		\caption{Convergence for increasing $L$.}\label{fig:Ishidirectconv}
	\end{subfigure}%
	\begin{subfigure}{0.52\textwidth}
		\includegraphics[height=0.195\textheight]{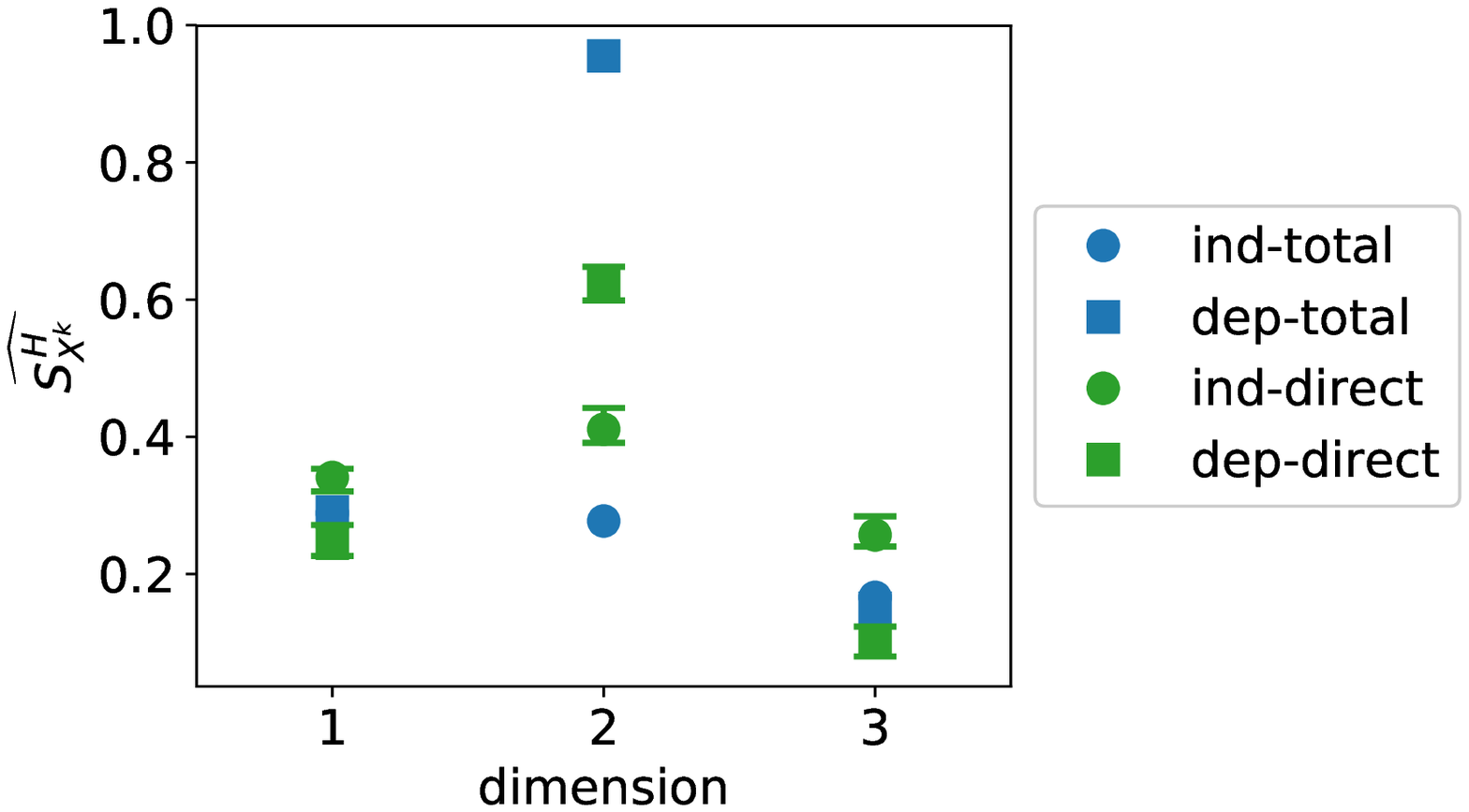}
		\caption{Comparison with the total indices.}\label{fig:Ishidirectcompare}	
	\end{subfigure}
	\caption{Direct indices for the Ishigami test case when computed with KDE.}\label{fig:Ishi5direct}
\end{figure}

We start with the KDE estimates in Figure \ref{fig:Ishi5direct}. On the left, we see that the estimates are relatively stable for increasing $L$ and the spread of the estimates decreases. On the right, we compare the estimates for $N=10^3$ for the direct indices to the estimates with $N=10^5$ for the total indices. We do not compute a reference value for the direct indices because of computational cost. For the dependent data, the indices work as expected: the total indices are larger than the direct \sis{}. For the independent data, this is not the case, as for variable 2 and 3 the direct \si{} is larger than the total \si{}. It is not immediately clear to us why this is the case, because for independent data, the total and direct \sis{} should give the same results.

\begin{figure}[ht!]
	\centering
	\begin{subfigure}{0.48\textwidth}
		\includegraphics[height=0.195\textheight]{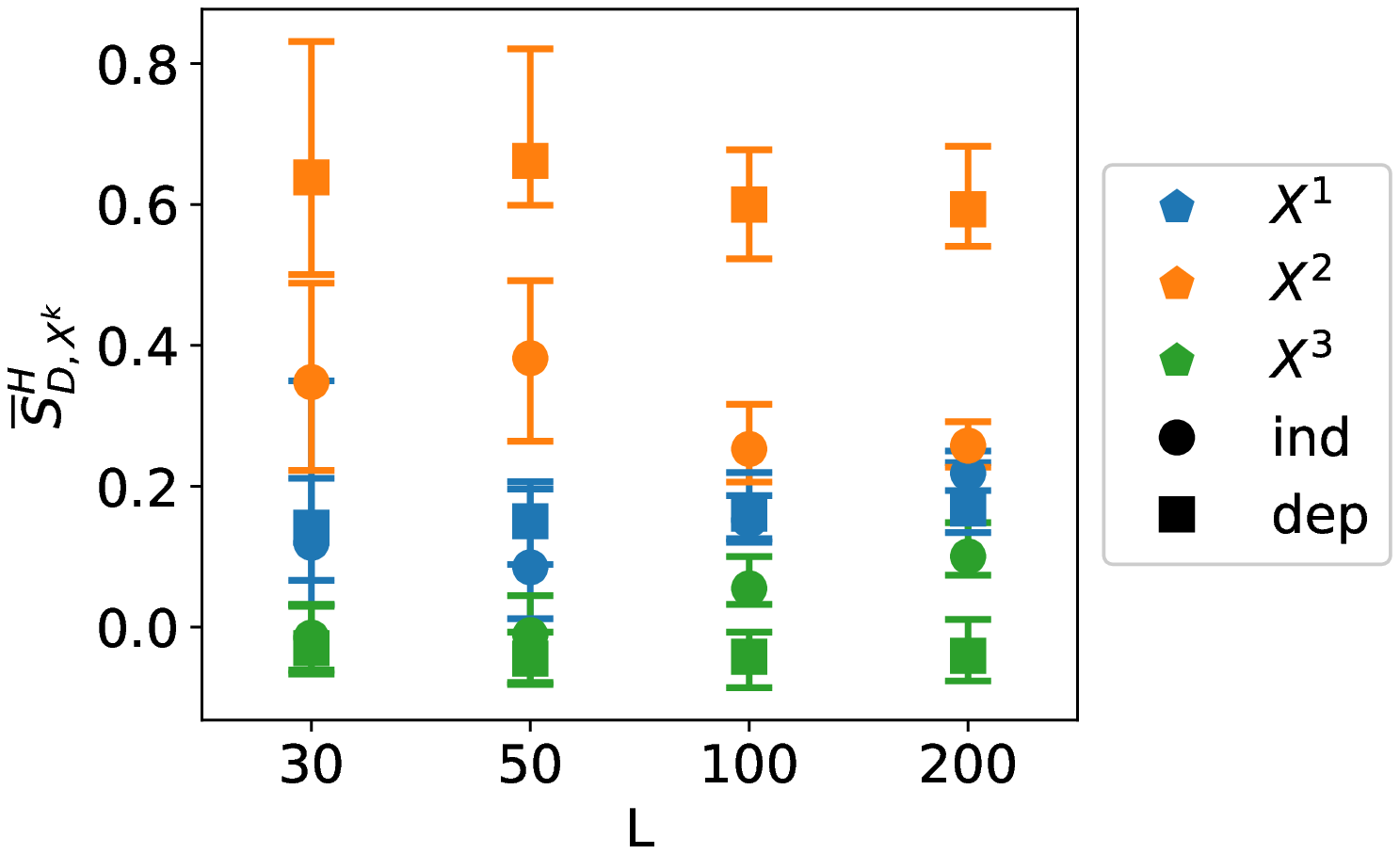}
		\caption{Convergence for increasing $L$.}\label{fig:IshidirectconvMST}
	\end{subfigure}%
	\begin{subfigure}{0.52\textwidth}
		\includegraphics[height=0.195\textheight]{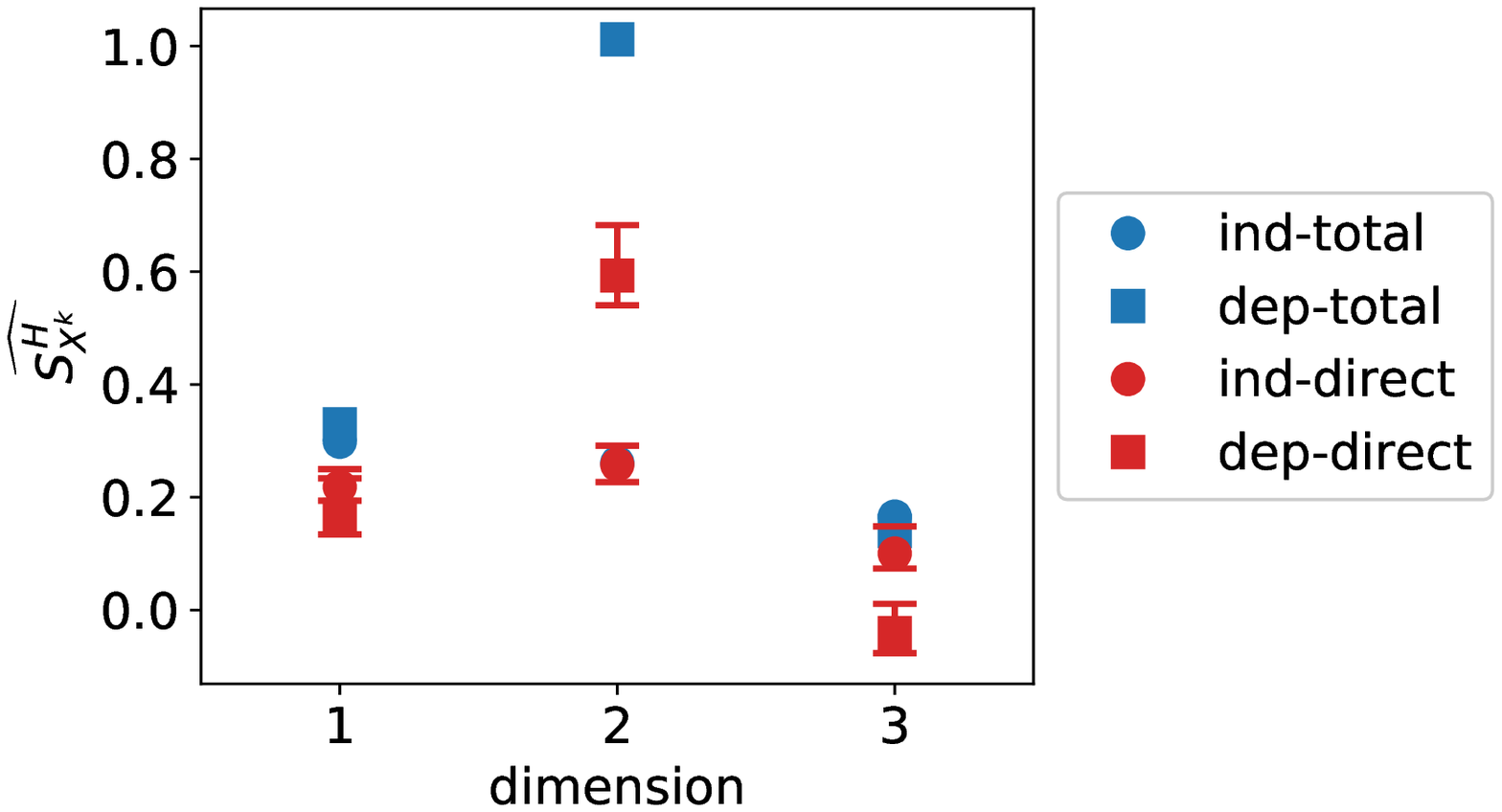}
		\caption{Comparison with the total indices.}\label{fig:IshidirectcompareMST}	
	\end{subfigure}
	\caption{Direct indices for the Ishigami test case when computed with MST.}\label{fig:Ishi5directMST}
\end{figure}

Figure \ref{fig:Ishi5directMST} shows on the left similar results for the stability and spread in the estimates as before with KDE. On the right, we see the total \sis{} are larger or equal than their direct \sis{} counterparts. For the independent data, the differences between the the direct and total \sis{} are small for variable 1 and 3 and nearly invisible for variable 2. Theoretically, this difference should be (numerically) zero. For the dependent data, we see the difference between direct and total \si{} is largest for variable 2, while variables 1 and 3 show a small difference. This is due to variable 2 being stronger dependent with the other two variables. 

\section{Conclusion}\label{sec:conclusion}
We proposed to use \gps{} in order to improve the estimates of divergence-based \sis{}. This is advantageous in cases where the number of available input-output samples is small, for example if the computational cost of each model evaluation needed to compute the output is high.

We compared the use of \gps{} to the well-established method of stochastic collocation combined with Lagrange interpolation. This method has several disadvantages in practice and is outperformed by the \gp{}-based methods in our experiments. The use of \gps{} also allowed us to propose (i) a new estimation method and (ii) a new type of \sis{}. This new estimation method for divergence-based \sis{} is based on minimum spanning trees and can be used in case the divergence used is the Hellinger distance. This estimation method has been used before to compute entropies and is numerically fast. The new type of \si{}, named direct \si{}, is especially useful when the input data is dependent. 

\section*{Acknowledgments}
This research is part of the EUROS programme, which is supported by NWO domain Applied and Engineering Sciences under grant number 14185 and partly funded by the Ministry of Economic Affairs.

\bibliographystyle{elsarticle-num}
\bibliography{references}
\end{document}